**Selection Rules and Channel Structure in a Base–Octave Model of Collatz Dynamics**

Katharina Lodders
Saint Louis, Missouri, USA.

Abstract

The Collatz iteration is governed by two distinct update rules, depending on the parity of the current iterate: $n_{i+1} = 3n_i + 1$ for odd $n_i$ and $n_{i+1} = n_i/2$ for even $n_i$. We show that these rules can be written equivalently as a single parity-controlled transformation, $n_{i+1}=((2s_i+1)(2k_i+s_i)+s_i)/2$, (where $n_i=2k_i+s_i$ and $s_i \in\{0,1\}$ is the parity of $n_i$.), yielding a uniform, step-aligned dynamical system in which parity variables are tracked explicitly. This reformulation removes the asymmetry of the traditional presentation and exposes structural regularities that are obscured when odd and even updates are treated separately.

Building on this unified rule, we introduce a base–octave decomposition, representing every integer uniquely in the form $n=B + 8(A-1)$ with $B \in\{1, \ldots,8\}$. The resulting dynamics separate into parity-dependent base transitions and affine updates of the octave index, inducing a finite directed transition skeleton lifted across scale levels. Refining the parity description yields a finite 128-state symbolic system that encodes all admissible transitions, including carry effects arising from higher-order parity inheritance.

Within this framework, we identify growth-permitting and decay-forcing channels and show that the only persistence mechanism—base-7 transitions in even octaves—is necessarily bounded by the 2-adic valuation of the octave index. An exhaustive enumeration of admissible return paths between persistence episodes establishes a non-positive drift in a logarithmic octave coordinate. Because of these finite-state constraints, trajectories are eventually confined to a contractive subnetwork associated with the terminal $1 \leftrightarrow 2$ cycle. The approach emphasizes structural organization and return-map methods, and provides a symbolic framework for analyzing parity-driven integer recurrences.

## 1. Introduction

The Collatz conjecture, introduced by Lothar Collatz in 1937, remains one of the most well-known unsolved problems in elementary number theory. It concerns the behavior of sequences generated from a positive integer $n \in \mathbb{N} = \{1,2,\ldots\}$ under the iterative rule

$$n_{i+1} = \begin{cases} n_i/2, & \text{if } n_i \text{ is even,} \\ 3n_i + 1, & \text{if } n_i \text{ is odd.} \end{cases}$$

The conjecture asserts that, regardless of the initial value $n > 1$, every such sequence eventually reaches the value 1. Despite extensive numerical verification and sustained theoretical effort, no general proof or disproof is currently known. The problem has been widely surveyed in the literature, and comprehensive historical accounts, references, and illustrative examples may be found in standard sources [1–2]. These will not be repeated here. Instead, the present work focuses on a structural reformulation of the iteration rule that enables a uniform accounting of sequence evolution and facilitates comparison across different starting values.

For reference, we briefly recall the traditional construction of Collatz sequences. Given an initial value $u \in \mathbb{N}$, one obtains a sequence $\{u_i\}_{i\geq 1} = (u_1, u_2, u_3, \ldots)$ by applying the above parity-dependent-rules at each step. When $u_i$ is odd, the transformation $3u_i + 1$ necessarily produces an even value, which is then followed by division by 2 in the subsequent step. Thus, for an odd iterate $u_i$, the combined transformation over two steps yields

$$u_{i+2} = \frac{3u_i + 1}{2}.$$

The parity of $u_{i+2}$ may again be odd or even, and similarly, the parity of iterates obtained from even starting values after division by 2 is not predictable a priori. Consequently, local parity information alone does not determine the subsequent evolution of a sequence, complicating attempts to analyze or compare Collatz trajectories in a systematic manner.

Empirically, Collatz sequences exhibit a general tendency toward convergence, but with substantial variation in their total stopping times. Some starting values give rise to unusually long trajectories, characterized by pronounced growth excursions and declines over relatively narrow step intervals. Well-known examples include sequences initiated at n =27 and 31, which require substantially more steps to reach 1 than neighboring initial values. If such excursions were unbounded, convergence would fail; understanding the structural mechanisms that permit large excursions and the constraints that limit them is therefore central to the analysis.

In what follows, a generalized formulation of the Collatz iteration is introduced that treats parity explicitly within a single update rule and applies a uniform step count across all sequences. This approach reveals regularities in the organization of trajectories, clarifies the origin of long and short stopping times, and provides a structured framework for describing the evolution of Collatz sequences from arbitrary initial values.

Accordingly, the aim of the paper is twofold: (i) to derive a finite, parity-controlled transition framework (selection rules and an explicit codebook) that organizes every trajectory as a path in a lifted finite-state system, and (ii) to establish a return-map inequality for the log–octave coordinate between successive base-7 persistence episodes. The return-map inequality is established in section 9 via a 128-state refinement and an exhaustive return-path 2-adic budget bound.

## 2. Generalized Form of Collatz Rules

Unlike earlier reformulations that track parity or valuation locally, the present framework yields a finite symbolic transition system governing all admissible base routing, with bounded growth certified via return-map drift rather than average estimates. The traditional Collatz iteration employs two distinct rules, $3n + 1$ for odd integers and $n/2$ for even integers. For analytical purposes, it is useful to place these operations on a common footing. This can be achieved by expressing both rules within a single parity-controlled transformation that applies uniformly to all integers.

Every integer $u \in \mathbb{N}$ can be written uniquely as

$$u_i = 2k_i + s_i, \qquad s_i \in \{0,1\}, \qquad (1)$$

where $s_i = 0$ corresponds to even integers and $s_i = 1$ to odd integers. The parameter $s_i$ thus encodes the parity of $u_i$.

Substituting (1) into the odd-update rule and making explicit the role of parity in the coefficients of the expression $3n + 1$, the transformation may be written as

$$u_{i+1} = (2s_i + 1)(2k_i + s_i) + s_i. \qquad (2)$$

For $s_i = 1$ (odd $u_i$), equation (2) reduces to

$$u_{i+1} = 3(2k_i + 1) + 1 = 3u_i + 1, \qquad (3)$$

which correctly reproduces the standard odd-update and yields an even integer. This value is then followed, in the traditional formulation, by a division by 2 in the subsequent step,

$$u_{i+2} = \frac{u_{i+1}}{2} = \frac{3u_i + 1}{2}.$$

The parity of $u_{i+1}$ is not determined by the parity of $u_i$ (but see below), and hence the subsequent evolution cannot be predicted locally.

For $s_i = 0$ (even $u_i$), equation (2) gives

$$u_{i+1} = 2k_i = u_i, \qquad (4)$$

so that the parity-selective factor leaves even integers unchanged under the $3n+1$ component. In the traditional Collatz sequence, the operative step for even integers is the division by 2,

$$u_{i+1} = \frac{u_i}{2},$$

Equation (2) makes explicit that the factor $2s_i+1$ performs a parity-selective multiplication: it equals 3 for odd integers and 1 for even integers. The additive term $+s_i$ ensures the parity change required for odd values while leaving even values unaffected. This formulation reveals that the traditional presentation implicitly treats odd and even integers asymmetrically in its step accounting. For odd integers, the operations $3n+1$ and the subsequent division by 2 are recorded as two distinct steps, whereas for even integers only the division by 2 is counted. As a result, the true operational steps across different sequences are not aligned.

This mismatch in step counting contributes to the apparent irregularities observed when Collatz sequences for consecutive starting values are compared at fixed step numbers. These discrepancies grow with increasing iteration count, but periodic structure emerges when sequences are examined under a uniform operational framework, as discussed in Section 3. These regularities motivate the introduction of a combined update rule.

Accordingly, the Collatz iteration may be written as a single parity-controlled transformation,

$$h_{i+1} = \frac{((2s_i + 1)(2k_i + s_i) + s_i)}{2}, \qquad h_i = 2k_i + s_i, \qquad (5)$$

where $h_i, k_i \in \mathbb{N}$, $s_i \in \{0,1\}$, and $i = 1,2,3, \ldots$ indexes the iteration steps. This generalized rule incorporates both the odd and even-updates into one expression and assigns a uniform step count to all transformations.

This form algins with the "accelerated" (or "short") Collatz map T(n)=n/2 if $n$ is even, and T(n)=(3$n$+1)/2 if $n$ is odd; this is the two-step compression of the standard map in which each odd step $n \mapsto 3n+1$ is followed immediately by one division by 2. This accelerated form follows from equation (5) and is used here throughout.

By eliminating the asymmetric treatment of odd and even integers, the generalized rule removes ambiguities in step accounting and allows Collatz sequences to be compared consistently across different initial values. This single update rule represents the accelerated Collatz iteration and is equivalent to two successive steps of the classical map when the iterate is odd. This uniform formulation exposes structural regularities that are obscured in the traditional representation and provides a foundation for deriving explicit parity-based constraints on admissible transitions in the subsequent sections. These constraints support a Lyapunov-style organization of ascent and descent mechanisms. The return-map analysis in section 9.6 supplies the final drift estimate.

The essential modification is to replace the traditional odd and even-updates by a single combined parity-sensitive operation as embodied in the generalized Collatz rule (equation (5)). When all operations are counted explicitly — even those that are operationally redundant (equation 4) — the traditional sequences expand, whereas using the combined single-step formulation of equation (5) yields the sequences shown in Table 1. Representative sequences under the traditional, fully expanded, and combined formulations are compared in the example $n = 20$; extended listings are omitted for brevity. In Table 1, each step corresponds to the same number of elementary operations applied to the original starting value. This uniform step counting exposes parity-dependent periodic patterns in the initial steps across all sequences, as discussed in Section 4.

Adopting the combined rule does not eliminate unusually long trajectories—such as those initiated at 27 and 31—but it alters their overall profiles by aligning step accounting and reducing artificial inflation in recorded step counts and intermediate values. The resulting sequences differ in their apparent heights and shapes while representing the same underlying integer dynamics. The mechanisms responsible for these outliers are analyzed below.

**Example**: As an illustrative example, the traditional Collatz sequence starting at $\underline{u}_1 = 20$ is

$$U_T(20) = \{u_i\}_{i=1}^{8}, (u_1, \dots, u_8) = (20,10,5,16,8,4,2,1).$$

which contains eight entries including the starting value. When every operation is recorded explicitly, the corresponding long table sequence becomes

$$U_L(20) = \{u_i\}_{i=1}^{14}, (u_1, \dots, u_{14}) = (20,20,10,10,5,16,16,8,8,4,4,2,2,1).$$

containing fourteen entries. In contrast, the combined step formulation (equation 5) yields

$$U_S(20) = \{u_i\}_{i=1}^{7}, (u_1, \dots, u_7) = (20,10,5,8,4,2,1).$$

with seven entries including the initial value (Table 1). In this representation, intermediate values such as 16 are absorbed into the combined (3u +1)/2 operation for odd inputs, resulting in a more compact and structurally transparent description of the sequence.

**Table1.** Examples of Collatz sequences under the generalized single-step formulation (equation (5)), showing uniform step counting and the resulting parity-dependent- periodic patterns across consecutive starting values. Entries marked "—" indicate values omitted after the sequence enters the terminal basin (the 1 ↔ 2 cycle).

| step n | 1 | 2 | 3 | 4 | 5 | 6 | 7 | 8 | 9 | 10 | 11 | 12 | 13 | 14 | 15 | 16 | 17 | 18 | 19 |
|---|---|---|---|---|---|---|---|---|---|---|---|---|---|---|---|---|---|---|---|
| **1** | 1 | 2 | 3 | 4 | 5 | 6 | 7 | 8 | 9 | 10 | 11 | 12 | 13 | 14 | 15 | 16 | 17 | 18 | 19 |
| **2** | — | 1 | 5 | 2 | 8 | 3 | 11 | 4 | 14 | 5 | 17 | 6 | 20 | 7 | 23 | 8 | 26 | 9 | 29 |
| **3** | — | — | 8 | 1 | 4 | 5 | 17 | 2 | 7 | 8 | 26 | 3 | 10 | 11 | 35 | 4 | 13 | 14 | 44 |
| **4** | — | — | 4 | — | 2 | 8 | 26 | 1 | 11 | 4 | 13 | 5 | 5 | 17 | 53 | 2 | 20 | 7 | 22 |
| **5** | — | — | 2 | — | 1 | 4 | 13 | — | 17 | 2 | 20 | 8 | 8 | 26 | 80 | 1 | 10 | 11 | 11 |
| **6** | — | — | 1 | — | — | 2 | 20 | — | 26 | 1 | 10 | 4 | 4 | 13 | 40 | — | 5 | 17 | 17 |
| **7** | — | — | — | — | — | 1 | 10 | — | 13 | — | 5 | 2 | 2 | 20 | 20 | — | 8 | 26 | 26 |
| **8** | — | — | — | — | — | — | 5 | — | 20 | — | 8 | 1 | 1 | 10 | 10 | — | 4 | 13 | 13 |
| **9** | — | — | — | — | — | — | 8 | — | 10 | — | 4 | — | — | 5 | 5 | — | 2 | 20 | 20 |
| **10** | — | — | — | — | — | — | 4 | — | 5 | — | 2 | — | — | 8 | 8 | — | 1 | 10 | 10 |
| **11** | — | — | — | — | — | — | 2 | — | 8 | — | 1 | — | — | 4 | 4 | — | — | 5 | 5 |
| **12** | — | — | — | — | — | — | 1 | — | 4 | — | — | — | — | 2 | 2 | — | — | 8 | 8 |
| **13** | — | — | — | — | — | — | — | — | 2 | — | — | — | — | 1 | 1 | — | — | 4 | 4 |
| **14** | — | — | — | — | — | — | — | — | 1 | — | — | — | — | — | — | — | — | 2 | 2 |
| **15** | — | — | — | — | — | — | — | — | — | — | — | — | — | — | — | — | — | 1 | 1 |

## 3. Net Changes

This section derives necessary bookkeeping identities that constrain the cumulative balance of increases and decreases along any Collatz trajectory. These identities describe what must hold if a trajectory converges, but do not by themselves establish convergence or rule out cycles. The structural machinery required to exclude nontrivial cycles is developed in Sections 5–9.

Let $\{H\}$ be a Collatz sequence generated by the generalized parity-controlled-rule in equation (5). If the sequence converges to 1, then the cumulative reductions arising from divisions by 2 must offset the cumulative increases generated by the $3n+1$ component acting on odd values. To quantify this balance, we consider the net change between successive steps.

The increment at step $i$ is

$$h_{i+1} - h_i = \frac{(2s_i+1)(2k_i+s_i)+s_i}{2} - (2k_i + s_i) \qquad (6)$$

## Even Steps

For even values ($s_i$ =0), equation (6) reduces to

$$h_{i+1} - h_i = -k_i, \qquad s_i = 0, \qquad (7)$$

which corresponds to the decrease produced by division by 2. Summing these contributions over the sequence yields the total decrement arising from even predecessors,

$$\Delta_{\text{even}} = \sum_{i=1}^{m-1} (-k_i)(1 - s_i) \qquad (8)$$

where the factor (1− $s_i$) ensures that only even terms are counted.

## Odd steps

For odd values ($s_i$ =1), equation (6) gives

$$h_{i+1} - h_i = k_i + s_i = k_i + 1, \qquad s_i = 1. \qquad (9)$$

The cumulative contribution from odd predecessors is therefore

$$\Delta_{\text{odd}} = \sum_{i=1}^{m-1} (k_i + s_i)s_i = \sum_{i=1}^{m-1} (k_i s_i + s_i^2). \qquad (10)$$

## Total net change

Combining even and odd contributions, the total net change over the m−1 operational steps is

$$
\begin{aligned}
\Delta_{\text{all}} &= \Delta_{\text{even}} + \Delta_{\text{odd}} \\
&= \sum_{i=1}^{m-1}(h_{i+1} - h_i) = \sum_{i=1}^{m-1}(2\,k_i s_i - k_i + s_i^2) \\
&= -\sum_{i=1}^{m-1} k_i + \sum_{i=1}^{m-1}(2\,k_i + s_i)s_i.
\end{aligned} \tag{11}
$$

The first summation term represents the cumulative reduction from divisions by 2, which apply to both odd and even values under the generalized rule. The second summation term represents the cumulative increase contributed exclusively by odd values, because the $3n+1$ component acts only when $s_i=1$.

Independently, the total change from the initial value $h_1$ to the terminal value $h_m$ satisfies

$$
h_m - h_1 = \sum_{i=1}^{m-1}(\,h_{i+1} - h_i) = \Delta_{\text{all}}. \tag{12}
$$

Substituting equation (11) into (12) yields

$$
h_m = h_1 - \sum_{i=1}^{m-1} k_i + \sum_{i=1}^{m-1}(2\,k_i + s_i)s_i = h_1 - \sum_{i=1}^{m-1} k_i + \sum_{i=1}^{m-1} h_i\,s_i. \tag{13}
$$

Equation (13) expresses the terminal value as the initial value reduced by the cumulative halving contributions and increased by the cumulative contributions from odd iterates. This relation describes the global balance of increases and decreases along a trajectory, rather than the behavior at individual steps.

## 3.1 Case Analysis of the Net Change

We now examine the possible outcomes implied by equation (13), which relates the initial and terminal values of a Collatz sequence through the cumulative balance of increases and decreases. Of particular interest is the special case relevant to the Collatz conjecture, namely $h_m = 1$.

### 3.1.1 Case 1: Steady-state, zero net change

Suppose that the cumulative increases and decreases exactly balance, so that

$$
h_m - h_1 = 0,
$$

or equivalently,

$$0 = \sum_{i=1}^{m-1} (-k_i + (2k_i + s_i)s_i). \qquad (14)$$

If this condition holds, then there exists a positive integer $h_i \in \mathbb{N}$ for which the generalized Collatz operation produces no net change over a finite interval of steps, implying $h_m = h_i$. Such a situation would correspond to the presence of a closed cycle within the sequence.

If, however, every Collatz sequence starting from any initial value $h_1 > 2$ reaches 1 in a finite number of steps, then no value $h_i > 2$ can occur more than once within a single sequence. The restriction $h_i > 2$ is essential: once the value 1 is reached, the sequence enters the trivial oscillatory cycle 1,2,1,2,..., allowing repeated occurrences of $h_i = 1$ and $h_i = 2$, hence the stopping time for sequence end is defined at the first occurrence of $h_m = 1$. The case $h_i = 0$ is excluded, as $0 \notin \mathbb{N}$ and would trivially satisfy equation (14).

The absence of nontrivial steady-state solutions is thus a necessary condition for convergence; showing that no such solutions exist requires the structural analysis of Sections 5–9.

### 3.1.2 Case 2: Net change by one

Consider next the case

$$h_m - h_1 = 1,$$

that is,

$$1 = \sum_{i=1}^{m-1} (-k_i + (2k_i + s_i)s_i). \qquad (15)$$

A net increase of one may occur locally within a sequence over several steps and arises, for example, when $h_i = 1$ and $h_m = 2$, which marks entry into the terminal oscillation 1,2,1,2,... . Because the Collatz iteration concerns initial values $h_1 > 1$, this case lies outside the boundary conditions of interest.

### 3.1.3 Case 3: Convergence to Unity ($h_m = 1$)

The case of primary interest for the Collatz conjecture is convergence to 1, that is,

$$h_m = 1.$$

Substituting this condition into equation (13) yields

$$1 = h_1 - \sum_{i=1}^{m-1} k_i + \sum_{i=1}^{m-1} (2\,k_i + s_i)s_i. \qquad (16)$$

Rearranging, this condition may be written equivalently as

$$h_1 - 1 = \sum_{i=1}^{m-1} k_i - \sum_{i=1}^{m-1} (2\,k_i + s_i)s_i. \qquad (17)$$

**Remark:** (Scope of this condition). Equation (17) is a necessary condition that any convergent trajectory must satisfy. It does not by itself exclude nontrivial cycles, which would require showing that no trajectory can return to a previously visited value. The structural machinery for ruling out cycles is developed in Sections 5–9.

Equation (17) expresses a necessary global condition for convergence to 1. It shows that the entire initial magnitude $h_1$, reduced by one, must be eliminated through the cumulative effect of the sequence operations. Thus, convergence requires not merely a balance between growth and decline, but a strict net reduction equal to the full initial value $h_1 - 1$.

The first summation on the right-hand side of (17) represents the cumulative reductions arising from all divisions by 2 applied throughout the sequence, regardless of parity. The second summation represents the cumulative growth contributed exclusively by odd iterates, when the $3n+1$ component acts for $s_i = 1$. Convergence to 1 therefore requires that the total halving effect exceed the total growth effect by exactly $h_1 - 1$.

In particular, it is not sufficient for growth and decline along a trajectory to balance on average. Any finite surplus of growth would leave the initial scale unreduced. Convergence requires a sustained predominance of halving operations, ensuring progressive reduction of the octave coordinate. Within the framework developed here, this requirement excludes steady states and non-trivial cycles as admissible terminal behaviors for trajectories with $h_1 > 1$. The return-map analysis of section 9.6 shows that the admissibility constraints enforce this dominance globally, yielding eventual confinement to the contractive regime and, in that case, possible convergence of all admissible trajectories.

Because the parity pattern $\{s_i\}$ and the intermediate values $h_i$ depend on the entire preceding path, neither parity nor monotonicity is known in advance. As a result, convergence cannot be established by induction or by local monotonicity arguments. Indeed, many sequences are highly non-monotonic and exhibit large but finite excursions over relatively short intervals before descending. Equation (17) provides a necessary global constraint: Every sequence that converges to 1 must satisfy this exact reduction condition. The net-change condition derived here is a necessary but not sufficient condition for convergence; additional structural constraints on admissible trajectories are required to rule out nontrivial cycles and are developed in the subsequent sections.

When combined with the parity-dependent selection rules governing admissible transitions derived below, this formulation together with a base and octave decomposition constrains the possible trajectories of Collatz sequences. The purpose of the base–octave decomposition (section 4) is not to simplify individual trajectories, but to expose a finite transition skeleton on which all trajectories are lifted.

Repeated growth excursions must be offset by sufficiently many halving steps so that the cumulative reduction removes not only the transient growth contributions but also the initial value to 1. The results of Sections 5–9 provide a finite-state structural framework for organizing and bounding these excursions. The return-map analysis in section 9.6 closes the remaining gap by establishing non-positive return drift.

## 4. Periodic Patterns Across Sequences with Consecutive Starting Values

When Collatz sequences are generated using the generalized parity-controlled rule (equation (5)), each iteration step represents the same number of elementary Collatz operations (cf. Table 1). This uniform step-accounting reveals systematic periodic structure in the parity patterns of iterates across sequences with consecutive starting values. The effect is most apparent in the initial iterations—particularly at steps $i$ = 2–4—but remains encoded, though increasingly diffuse, at later steps.

More specifically, when sequences are ordered by their initial values $h_1(n)$=1,2,..., the parities of the iterates $h_2(n)$ through $h_4(n)$ display a regular periodic structure (Appendix B gives examples). This structure extends across two consecutive octaves of starting values and reflects the explicit parity dependence of the Collatz update rule. At this stage, these correlations describe systematic parity patterns rather than exact numerical outcomes; however, they form the basis for the parity-dependent selection rules derived later, which permit constrained prediction of subsequent base transitions and iterate behavior.

Although these periodic patterns become progressively less apparent at higher step numbers, they remain detectable through at least step $i$ = 5. The underlying reason is that all sequences whose starting values differ by eight share identical parity patterns at corresponding steps. Because $8 = 2^3$, this periodicity reflects the three-bit parity structure inherent in the binary decomposition of integers and the parity-controlled nature of the update rule. Through equation (5), these shared parity patterns induce structured and repeatable variations in subsequent iterates, which are later formalized as explicit selection rules.

Thus, the structural analysis of Collatz sequences may be reduced from infinitely many starting values to a finite set of eight base cases. These base cases correspond to the residue classes modulo 8 and span two consecutive octaves of starting values. The parity-dependent selection rules governing transitions between base classes are derived from this finite set in the subsequent analysis.

To formalize this decomposition, each integer $h_i$ is written uniquely in the form

$$h_i = B_i + 8(A_i - 1), \qquad B_i \in \{1, \dots, 8\},\ A_i \in \mathbb{N}, \qquad (18)$$

or equivalently,

$$h_i = 8A_i - (8 - B_i), h_i > 0. \qquad (19)$$

Here $B_i$ denotes the base value identifying one of the eight base sequences, and $A_i$ denotes the corresponding octave index. This representation assigns every integer to exactly one base value $B_i$ and one octave index $A_i$.

The base value associated with a given $h_i$ is

$$B_i(h_i) = ((h_i - 1) \, mod \, 8) + 1, \qquad (20)$$

where the shift by one avoids the undesired zero residue. The corresponding octave index is

$$A_i(h_i) = 1 + \frac{h_i - B_i}{8}. \qquad (21)$$

The eight base values $B = 1,\dots,8$ define eight vertical base sequences, each containing exactly one element in every octave. The octaves themselves, indexed by $A = 1,2,\dots$, partition the natural numbers into consecutive blocks of eight. Because the term $8(A_i - 1)$ is always even, all integers belonging to the same base sequence share the same parity as their base value $B_i$. This parity invariance is a direct consequence of equation (18) and is fundamental to the parity-dependent transition structure developed in the subsequent analysis.

For example, the integers $h = 1,9,17, \dots$ all share the same base value $B = 1$ while occupying successive octaves $A = 1,2,3, \dots$. Each such set forms a vertical base sequence extending across all octaves. When referring to the eight base classes $B \in \{1, \dots, 8\}$, we use the labels B1, B2,.., B8 for these classes. Subscripts $i$ are reserved exclusively for iteration indices within Collatz sequences and are not used to label base classes.

To refine parity information beyond a single iteration step, it is useful to resolve not only the parity of the base value $B$ and the octave index $A$, but also the parities of their associated subcomponents. Accordingly, we decompose the base value as

$$B_i = 2k_{b,i} + s_{b,i}, \quad (22)$$

where $s_{b,i}$ equals the parity of $B_i$ because the term $8(A_i - 1)$ in the decomposition $h_i = B_i + 8(A_i - 1)$ is always even. Likewise, we decompose the octave index as

$$A_i = 2k_{a,i} + s_{a,i}, \quad (23)$$

where $s_{a,i}$ records the parity of the octave number $A_i$. The periodic structure observed in the early iteration steps is naturally organized into blocks of two octaves, and it is convenient to record the parity of the octave index explicitly. The variable $s_{a,i}$ therefore serves as a bookkeeping device that distinguishes between odd and even octaves and will be used in the subsequent derivation of the parity-dependent selection rules.

To capture finer structure within octaves, we further decompose the half-octave index $k_a$ (sometimes referred to as the quartet index) by

$$k_a = 2k_q + s_q, \quad s_q \in \{0,1\}. \qquad (24a)$$

And furthermore, for later analysis:

$$k_q = 2k_r + s_r, \quad s_r \in \{0,1\}. \qquad (24b)$$

We resolve the sub-parity of the base coefficient $k_b$ via

$$k_b = 2k_c + s_c, \quad s_c \in \{0,1\}. \qquad (25)$$

These refinements introduce no new dynamics; they serve to record additional binary information required to propagate parity deterministically under the generalized Collatz update.

Substituting the $A, B$ decomposition $h_i = B_i + 8(A_i - 1)$ into the generalized parity-controlled update rule (5) expresses the evolution of Collatz sequences entirely in terms of base values and octave indices. Each iteration produces a new pair $(B_{i+1}, A_{i+1})$, where the base transition is determined by the parity structure of the current state, but the octave index updates according to a linear rule associated with the selected transition. In this formulation, the sequence pathway is governed by finite parity data, whereas the magnitude of $A$ encodes scale information carried forward along that pathway.

$$h_i = B_i + 8(A_i - 1) = (2k_{b,i} + s_{b,i}) + 8(2k_{a,i} + s_{a,i} - 1). \qquad (26)$$

Because the parity of $h_i$ is determined solely by $B_i$, the parity variable in the generalized update (equation (5)) is $p_i = s_{b,i}$. Substituting the base–octave decomposition (equation (26)) into equation (5) therefore yields

$$\begin{aligned} h_{i+1} &= \frac{(2s_{b,i} + 1)(2k_{b,i} + s_{b,i} + 8(A_i - 1)) + s_{b,i}}{2} \\ &= \frac{(2s_{b,i} + 1)(2k_{b,i} + s_{b,i}) + s_{b,i}}{2} + \frac{(2s_{b,i} + 1)\, 8(2k_{a,i} + s_{a,i} - 1)}{2}. \end{aligned} \qquad (27)$$

Equations (18)–(21), together with (26)–(27), apply to all starting values $h_1$ and account for the periodic step value patterns observed empirically in the early iterates when the parity-controlled update is applied to consecutive initial values $h_1(n)$=1,2,..., (cf. Table 1). We return to the quantitative consequences of these decompositions when deriving the transition (selection) rules.

With the base–octave decomposition and its parity refinements fixed in equations (18)–(27), we now introduce compact notation for separating contributions by base parity and octave parity in subsequent derivations.

## Indicator and sign factors (notation only)

In later derivations it is convenient to separate contributions according to octave parity and base parity without repeatedly restating case distinctions. For this purpose we use the indicator

$$1 - s_{a,i}, \qquad (28a)$$

which selects even octaves ($s_{a,i} = 0$), and the sign factor

$$1 - 2s_{b,i} \in \{+1, -1\}, \qquad (28b)$$

which distinguishes even and odd base parity. These factors introduce no new dynamics; they serve only as shorthand for case distinctions already implicit in the base–octave parity variables defined above (cf. (22)–(26)).

### Mod8 structure and the "half-octave" effect

The second term in (27) is always a multiple of 4. Its contribution modulo 8 depends only on the parity of $A_i$, i.e. on $s_{a,i}$:

$$4(2s_{b,i}+1)(A_i-1) \equiv \begin{cases} 0(\text{mod}8), & s_{a,i}=1 \quad (A_i \text{ odd}) \\ 4(\text{mod}8), & s_{a,i}=0 \quad (A_i \text{ even}) \end{cases} \qquad (29)$$

because $(2s_{b,i}+1)$ is always odd. Equivalently,

$$4(2s_{b,i}+1)(A_i-1) \equiv 4(1-s_{a,i})(\text{mod}8). \qquad (30)$$

Thus, within each octave, the contribution of the octave term to $h_{i+1}$ modulo 8 depends only on whether the iterate lies in the lower or upper half of the octave.

### Defining $B_{i+1}$ and $A_{i+1}$

To obtain the new base value and octave index after one update, we must apply the *same* canonical mapping as in (20)–(21). In particular, $B_{i+1}$ is defined by

$$B_{i+1} = ((h_{i+1}-1)\, mod\, 8)+1, \qquad (31)$$

and then

$$A_{i+1} = 1+\frac{h_{i+1}-B_{i+1}}{8}. \qquad (32)$$

Equations (27)–(32) apply to all starting values $h_1$ and explain why iterating the parity-controlled rule over consecutive initial values $h_1(n)$ = 1,2, ..., produces periodic behavior in the early step values $h_2(n)$, $h_3(n)$,..., (cf. Table 1). In particular, the residue class $B_{i+1}$ depends on (i) the "base part" $B_i$ and (ii) the octave parity $s_{a,i}$ through the mod-8 shift identified in (30). Enumerating the possibilities for $B_i \in \{1, ...,8\}$ and $s_{a,I} \in \{0,1\}$ yields 8×2 = 16 parity-controlled transition cases ("selection rules") governing movement between base sequences from one step to the next.

## 5. Selection Rules

The permitted transitions ("selection rules") governing the evolution of the base–octave pair $(B_i, A_i)$ follow from the decomposition in equation (18) together with the generalized Collatz update (equation (5)). Because the term $8(A_i-1)$ is even, the parity of $h_i$ is determined entirely by the parity of the base value $B_i$; this choice isolates parity information in the base variable and simplifies the subsequent case analysis.

**Lemma 5.1** (Parity isolation in the base–octave decomposition)

Let an iterate $h_i$ be written uniquely in base–octave form as in (18). Then the parity of $h_i$ is exactly the parity of $B_i$.

**Proof.** By construction, the term $8(A_i - 1)$ is even. Hence

$$h_i \equiv B_i(\mathrm{mod}2),$$

so the parity of $h_i$ is determined entirely by the parity of the base value $B_i$. ■

## Lemma 5.2 (Even-octave modular shift).

Let $h_i = B_i + 8(A_i - 1)$ be written in base–octave form, and suppose the generalized Collatz update produces an expression of the form

$$h_{i+1} = f(B_i) + c\,(A_i - 1),$$

where *f(B_i)* depends only on $B_i$ and $c \in \{4,12\}$.
If $A_i$ is even, then $A_i - 1$ is odd and

$$c(A_i - 1) \equiv 4(\mathrm{mod}8).$$

Consequently, the canonical base value $B_{i+1}$ differs from the residue of *f(B_i)* by an additive shift of +4modulo 8.

**Proof.** If $A_i$ is even, then $A_i - 1$ is odd. For $c = 4$ or $c = 12$, one has

$$c(A_i - 1) \equiv 4(\mathrm{mod}8).$$

Under the canonical mapping (31)-(32), this contributes a $+4$ shift to the base residue. ■

By Lemma 5.1, parity information is encoded in the base variable $B_i$, while the octave index $A_i$ governs modular shifts without affecting parity. To preserve uniqueness of the base–octave representation after each update, the next base and octave are defined canonically by equations (31)-(32). The update therefore splits naturally into four cases according to the parity of $B_i$ and the parity of $A_i$.

### Case 1 (even base, odd octave): $B_i$ even, $A_i$ odd

By Lemma 5.1, $B_i$ even implies $h_i$ is even. The generalized Collatz update therefore reduces to halving, i.e. $h_{i+1}=h_i/2$.

$$h_{i+1} = h_i/2.$$

Using the base–octave decomposition (18),

$$h_{i+1} = \frac{B_i + 8(A_i - 1)}{2} = \frac{B_i}{2} + 4(A_i - 1). \qquad (33)$$

When $A_i$ is odd, the quantity $A_i - 1$ is even, and hence $8(A_i - 1)$ is a multiple of 16 (in particular, a multiple of 8). No modular shift of the base residue occurs. Under the canonical mapping (31)–(32),

$$B_{i+1} = \frac{B_i}{2}, A_{i+1} = 1 + \frac{A_i - 1}{2} \qquad (34)$$

The permitted base transitions for odd octaves are therefore

$$B2 \to B1, \quad B4 \to B2, \quad B6 \to B3, \quad B8 \to B4.$$

### Case 2 (even base, even octave): $B_i$ even, $A_i$ even

As in Case 1, Lemma 5.1 implies that hi is even and the halving form applies. If $A_i$ is even, then $A_i - 1$ is odd. By Lemma 5.2, the term $4(A_i - 1)$ contributes a residue of +4modulo 8. Thus, the base residue undergoes a +4 shift under the canonical mapping. Since $\frac{B_i}{2} \in \{1,2,3,4\}$, no wraparound occurs, and

$$B_{i+1} = \frac{B_i}{2} + 4, \; A_{i+1} = \frac{A_i}{2}. \qquad (35)$$

The permitted base transitions for even octaves are

$$B2 \to B5, \quad B4 \to B6, \quad B6 \to B7, \quad B8 \to B8.$$

### Case 3 (odd base, odd octave): $B_i$ odd, $A_i$ odd

By Lemma 5.1, $B_i$ odd implies $h_i$ is odd. The generalized Collatz update yields

$$h_{i+1} = \frac{3h_i + 1}{2} = \frac{3B_i + 1}{2} + 12(A_i - 1). \qquad (36)$$

When $A_i$ is odd, the quantity $A_i - 1$ is even, so $12(A_i - 1)$ is a multiple of 8 and produces no modular shift. Consequently, the next base value is determined solely by the residue of $(3B_i + 1)/2$ under the canonical mapping (31)-(32),

$$B_{i+1} = ((\frac{3B_i + 1}{2} - 1)\, mod\, 8) + 1. \qquad (37)$$

Evaluating for $B_i \in \{1,3,5,7\}$ gives

$$B1 \to B2, \qquad B3 \to B5, \qquad B5 \to B8, \qquad B7 \to B3.$$

The corresponding octave index $A_{i+1}$ is then fixed uniquely by the canonical definition (31)-(32).

## Case 4 (odd base, even octave): $B_i$ odd, $A_i$ even

Equation (36) again applies. If $A_i$ is even, then $A_i-1$ is odd. By Lemma 5.2, the octave term contributes a residue shift modulo 8. Hence the base residue obtained from the odd-step numerator is shifted under the canonical mapping.

$$B_{i+1} = \left(\left(\frac{3B_i + 1}{2} + 4 - 1\right) mod\, 8\right) + 1. \qquad (38)$$

Evaluating for $B_i \in \{1,3,5,7\}$ yields

$$B1 \to B6, \qquad B3 \to B1, \qquad B5 \to B4, \qquad B7 \to B7.$$

As in Case 3, the octave index $A_{i+1}$ is determined uniquely by (31)-(32).

## Summary: the 16 selection rules

Cases 1–4 exhaust all combinations of base parity and octave parity. Because $B_i$ has eight possible values and the octave parity has two possibilities, these rules generate 8×2 =16 parity-controlled-transition cases. Lemma 5.1 isolates parity in the base variable, whereas Lemma 5.2 accounts for the +4modular shifts induced by even octaves. Together, they fully determine the admissible evolution of the base–octave pair under the generalized Collatz update.

**Proposition 5.3 (Exhaustiveness of the 16 base selection rules).**

For every iterate $h_i$ written uniquely as $h_i=B_i+8(A_i-1)$, the next base value $B_{i+1}$ is determined uniquely by the pair $(B_i,s_{a,i})$, and the mapping $(B_i,s_{a,i})\mapsto B_{i+1}$ is exactly the 16-case selection-rule list derived in Cases 1–4.

**Proof.** The generalized update (5) is deterministic on integers. By Lemma 5.1, the parity input $p_i$ in (5) is determined by $B_i$. The only remaining discrete choice affecting the residue class modulo 8 is whether the octave index is odd or even, i.e. $s_{a,I} \in \{0,1\}$, because even-octave updates produce the fixed modulo-8 shift of Lemma 5.2 whereas odd-octave updates do not. Therefore, the next residue class $h_{i+1}$ (mod 8) is determined uniquely by $(B_i, s_{a,i})$. Applying the canonical remapping (31)–(32) then yields a unique base representative $B_{i+1} \in \{1,\dots,8\}$. Cases 1–4 enumerate all possibilities for $B_i$ (eight values) and $s_{a,i}$ (two values), hence exhaust all 16 parity-controlled base-transition cases. ■

**Remark 5.4 (From base selection rules to the exhaustive refined codebook in Appendix A1).** Proposition 5.3 fixes the directed base-transition skeleton. To lift this skeleton to a deterministic state update, one also tracks the octave index update. As shown in Section 7.2, for each of the 16 base-transition cases the octave coordinate updates by an affine rule $A_{i+1}=mA_i+c$ (with $m \in\{1,3\}$ and an offset $c$ determined by the case) followed by the canonical base–octave remapping (20)–(21). The only additional splitting needed at the level of the finite parity descriptor Σ arises when the affine numerator introduces a carry across a power-of-two boundary in the binary expansion of $A_i$, which can change the inheritance of $s_{a,i+1}$ and downstream sub-parity bits. For each base-transition case this carry/no-carry distinction yields finitely many admissible refined outcomes. Appendix A1 lists these admissible carry/branch outcomes explicitly; hence, once Propositions 5.3 and the affine-update forms of Section 7.2 are accepted, Table A1 is exhaustive for the refined transitions.

## 6. Selection Rules as a Directed Transition Graph

The parity-controlled update rule derived in Section 5 induces a finite set of permitted transitions between base classes. In the base–octave representation

$$h_i = B_i + 8(A_i - 1),$$

the next base value $B_{i+1}$ is determined by the current base value $B_i$ together with the parity of the octave index $A_i$, via the canonical mapping

$$B_{i+1} = ((h_{i+1} - 1)\, mod\, 8) + 1,$$

in combination with the four parity cases established in Section 5.

Because the octave parity ($s_{a,i}$) can take two values (odd or even), each base value $B_i$ admits two distinct outgoing base transitions: one corresponding to odd octaves and one corresponding to even octaves. Likewise, each base value typically has two incoming transitions, arising from the two possible octave parity regimes. The selection rules therefore define a directed transition structure on the finite set of base values $\{1, \dots, 8\}$.

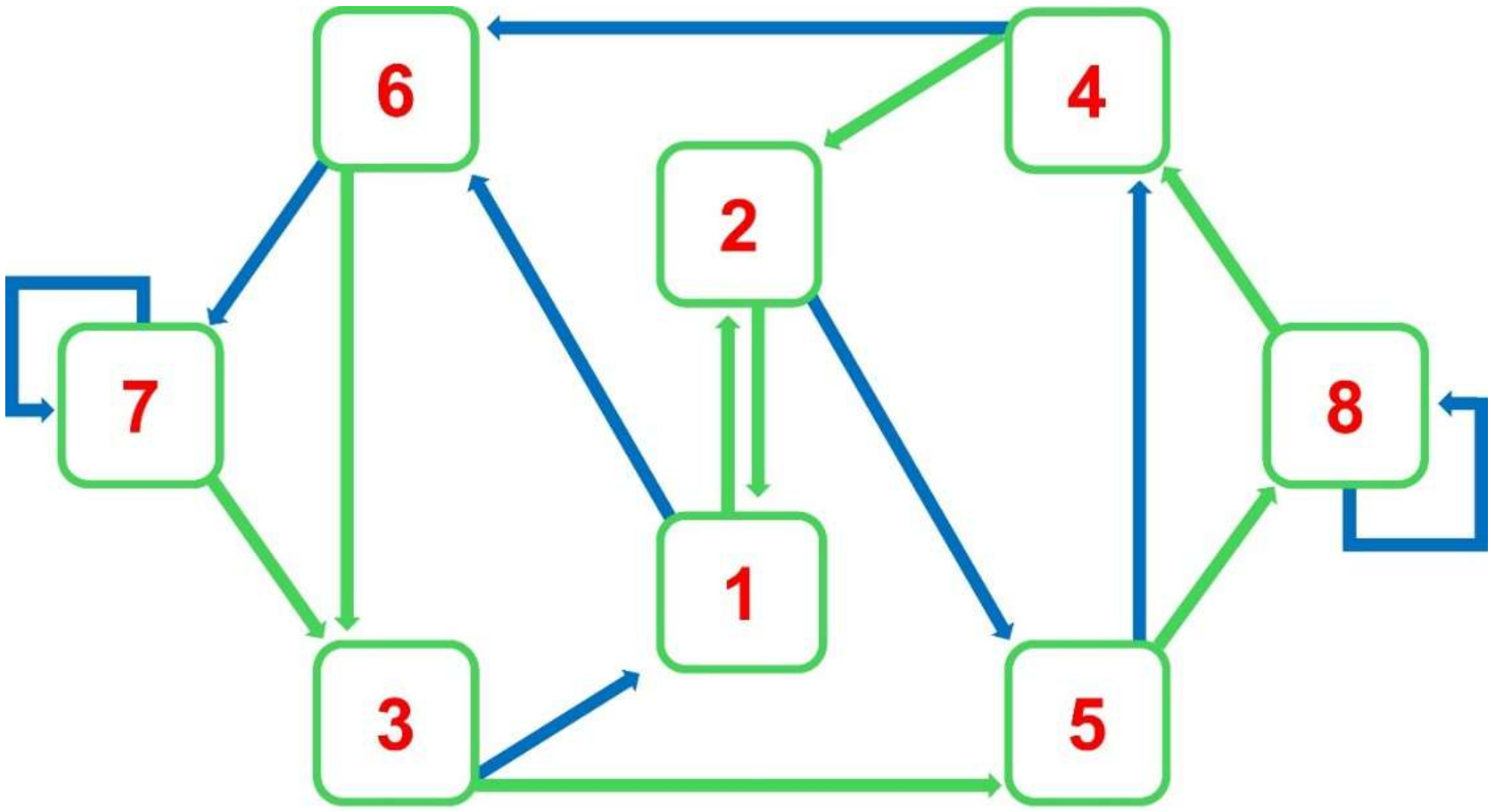

**Figure 1.** Permitted transitions induced by the Collatz rule when each iterate $h_i$ is uniquely expressed in base–octave form, $h_i = B_i + 8(A_i - 1)$, with base value $B_i \in \{1, \dots, 8\}$ and octave index $A_i \in \mathbb{N}$ (see text for examples). The Collatz update allows transitions only between specific base sequences, determined jointly by the parity of the base value $B_i$ and the parity of the octave $A_i$ containing it (cf. the selection rules in Section 5). Transitions originating from odd-numbered octaves (e.g., $h_i = 1 - 8, 17 - 24, \dots$) are shown in green, and transitions originating from even-numbered-octaves (e.g., $h_i = 9 - 16, 25 - 32, \dots$) are shown in blue. Although the allowed outgoing transitions differ between odd and even octaves, the same destination base *values* may be reached from both octave parities, possibly at different octave levels.

It is natural to represent this structure as a directed transition graph**,** in which vertices correspond to base values and directed edges represent admissible transitions $B_i \to B_{i+1}$. This representation is shown in Figure 1. The two parity regimes are distinguished by edge coloring: Permitted transitions originating from odd octaves are shown in green, and permitted transitions originating from even octaves in blue. The arrow directions encode the permitted base transitions derived in Section 5. The corresponding update of the octave index $A_i \to A_{i+1}$ is not shown explicitly in the graph but is determined uniquely by the canonical mapping once $h_{i+1}$ and $B_{i+1}$ are fixed.

The resulting directed graph provides a compact visualization of the 16 parity-controlled-base transition cases (one for each pair ($B_i$, $s_{a,i}$)). While the base transition graph is finite, the full dynamics occur on a lifted state space in which each base node is replicated across infinitely many octave layers. Figure 1 may therefore be viewed as the projection of a layered transition structure onto the base index alone.

Table 2 illustrates the action of the parity-controlled Collatz update on the eight base values across the first two octaves. The transitions confirm the analytically derived selection rules; no additional cases occur for higher octaves due to periodicity modulo 8.

**Table 2.** Base–Octave Decomposition and Generalized Collatz Update for $h_1 = 1$–$16$

| $A_i$ | $h_i$ | $B_i$ | $k_{b,i}$ | $s_{b,i}$ | Parity Change | $h_{i+1}$ | $B_{i+1}$ | $k_{b,i+1}$ | $s_{b,i+1}$ | $A_{i+1}$ | $s_{a,i+1}$ |
|---|---|---|---|---|---|---|---|---|---|---|---|
| $A$=1, $s_{a,1}$ = 1 | 1 | 1 | 0 | 1 | o → e | 2 | 2 | 1 | 0 | 1 | 1 |
| | 2 | 2 | 1 | 0 | e →o | 1 | 1 | 0 | 1 | 1 | 1 |
| | 3 | 3 | 1 | 1 | o →o | 5 | 5 | 2 | 1 | 1 | 1 |
| | 4 | 4 | 2 | 0 | e →e | 2 | 2 | 1 | 0 | 1 | 1 |
| | 5 | 5 | 2 | 1 | o → e | 8 | 8 | 4 | 0 | 1 | 1 |
| | 6 | 6 | 3 | 0 | e →o | 3 | 3 | 2 | 1 | 1 | 1 |
| | 7 | 7 | 3 | 1 | o →o | 11 | 3 | 5 | 1 | 2 | 0 |
| | 8 | 8 | 4 | 0 | e →e | 4 | 4 | 2 | 0 | 1 | 1 |
| $A$ = 2, $s_{a,2}$ = 0 | 9 | 1 | 0 | 1 | o → e | 14 | 6 | 7 | 0 | 1 | 1 |
| | 10 | 2 | 1 | 0 | e →o | 5 | 5 | 2 | 1 | 1 | 1 |
| | 11 | 3 | 1 | 1 | o →o | 17 | 1 | 8 | 1 | 3 | 1 |
| | 12 | 4 | 2 | 0 | e →e | 6 | 6 | 3 | 0 | 1 | 1 |
| | 13 | 5 | 2 | 1 | o → e | 20 | 4 | 10 | 0 | 3 | 1 |
| | 14 | 6 | 3 | 0 | e →o | 7 | 7 | 3 | 1 | 1 | 1 |
| | 15 | 7 | 3 | 1 | o →o | 23 | 7 | 3 | 1 | 2 | 0 |
| | 16 | 8 | 4 | 0 | e →e | 8 | 8 | 4 | 0 | 1 | 1 |

Note: The first two octaves ($A$=1, $A$=2) are shown for illustration; higher octaves follow the same update rules and differ only in the scale parameter $A$.

## Internal channels and special behavior of *B* = 7 and *B* = 8

The transition structure is parity-asymmetric: for most bases, the odd-octave outgoing edge and the even-octave outgoing edge lead to different destinations. Two base values exhibit an additional special feature:

$B$ = 7 has an internal even-octave-channel $B7 \to B7$, whereas in odd octaves $B7 \to B3$.

$B$ = 8 has an internal even-octave-channel $B8 \to B8$, whereas in odd octaves $B8 \to B4$.

These internal channels do not create new base classes, but they allow trajectories to remain in the same base class while changing octave level. This mechanism provides a natural structural explanation for why certain trajectories may spend extended periods in high-value regions (notably those frequently-visiting the $B$ = 7 class in even octaves), whereas others collapse quickly through repeated $B$ = 8 even-octave-transitions. (The quantitative effect on total stopping time is discussed later.)

### Counting destinations: 16 base transitions and 32 (*B*,*A*) destinations

There are 8 base values and 2 octave parities, hence 16 parity-controlled base-transition cases (one for each $(B, \mathrm{par}\,(A))$). However, the *destination octave parity* is not fixed solely by the base transition: depending on the higher-order components of $A_i$, the same base transition can land in either an odd or an even octave. Consequently, each of the 16 base transitions can split into two (*B*,*A*) destinations, giving up to 32 distinct $(B_{i+1}, \mathrm{par}\,(A_{i+1}))$ outcomes.

### Worked example: the trajectory from $h_1$ =7

As an illustrative example, consider the sequence starting at $h_1$ = 7. Using $h_i$ = $B_i$ +8($A_i$ −1), one obtains

$$\begin{aligned} 7\,(B = 7, A = 1) \quad &\to\ 11\,(B = 3, A = 2)\ \to\ 17\,(B = 1, A = 3)\ \to\ 26\,(B = 2, A = 4) \\ &\to\ 13\,(B = 5, A = 2)\ \to\ 20\,(B = 4, A = 3)\ \to\ 10\,(B = 2, A = 2) \\ &\to\ 5\,(B = 5, A = 1)\ \to\ 8\,(B = 8, A = 1)\ \to\ 4\,(B = 4, A = 1) \\ &\to\ 2\,(B = 2, A = 1)\ \to\ 1\,(B = 1, A = 1). \end{aligned}$$

The transitions 11 → 17, 26 → 13 → 20, and 10 → 5 originate from even octaves and follow the even-octave rules ($B3 \to B1$), ($B2 \to B5$), ($B5 \to B4$), and ($B2 \to B5$), respectively. In contrast, the transitions 17 → 26 and 20 → 10 originate from odd octaves and follow the odd-octave rules ($B1 \to B2$) and ($B4 \to B2$). This illustrates how the same base value can induce different base transitions depending on the parity of the octave in which it occurs.

## 7. Deterministic State Description and Finite Transition Structure

The selection rules derived in Section 5 and visualized as a directed base transition graph in Section 6 determine how the base index $B_i$ evolves from step to step. To obtain a complete and deterministic description of the Collatz iteration, however, additional information beyond the base value and octave index is required. This information is captured by a refined parity state that encodes higher order-parity data governing the evolution of the octave index.

**Terminology.** We use *trajectory* to mean the forward sequence $(h_i)_{i\geq 0}$ defined by $h_{i+1} = T(h_i)$. We use *pathway* to mean the induced sequence of discrete states (e.g., base indices $B_i$ or refined states $\Sigma_i$) along the finite transition graph/codebook. No cyclic behavior is implied by either term.

## 7.1 Refined state representation

At each iteration step $i$, we associate to the iterate $h_i$ the state vector

$$\Sigma_i = (B_i,\ s_{b,i},\ s_{c,i},\ s_{a,i},\ s_{q,i}),$$

where:

- $B_i \in \{1, \dots, 8\}$ is the base index in the decomposition $h_i = B_i + 8(A_i - 1)$;
- $s_{b,i} \in \{0,1\}$ is the parity of the base value $B_i$;
- $s_{a,i} \in \{0,1\}$ is the parity of the octave index $A_i$;
- $s_{c,i}$ and $s_{q,i}$ are parity indicators (equations 24, 25) arising from the quotient terms in the generalized Collatz update and encode higher-order binary information not contained in $s_{b,i}$ or $s_{a,i}$ alone.

The introduction of $s_{c,i}$ and $s_{q,i}$ allows the update of the octave parity and the base parity at the next step to be expressed deterministically.

## 7.2 Deterministic Update Structure and Conditional Parity Inheritance

We now clarify the extent to which the refined parity variables evolve deterministically under the generalized Collatz update. The full state of the iteration at step $i$ is the pair

$$(A_i, \Sigma_i),$$

where $A_i$ is the octave index and $\Sigma_i = (B_i, s_{b,i}, s_{c,i}, s_{a,i}, s_{q,i})$ is the refined parity state, where $B_i$ is the base index, $s_{b,i}$ its parity, $s_{a,i}$ the parity of the octave index $A_i$, and $s_{c,i}, s_{q,i}$ encode higher-order parity information arising from the quotient terms in the update.

The state $\Sigma_i$ is intended as a finite parity-augmented descriptor governing the transition skeleton on base classes, whereas the octave index $A_i$ is treated as a separate, unbounded level coordinate, consistent with the description of the dynamics as finite parity states coupled to an infinite family of octave levels.

Although $\Sigma_i$ records several parity features of the octave level, it does not, by itself, encode all higher-order binary information in $A_i$. In particular, the evolution of some next-step parity bits may depend on whether the affine octave update introduces a carry across a power-of-two boundary in the binary expansion of $A_i$. Accordingly, determinism holds for the full state pair $(A_i,\Sigma_i)$, while the projected finite descriptor $\Sigma_i$ may require an explicit carry/branch case to resolve certain parity inheritances.

The base parity update satisfies the exact identity

$$s_{b,i+1} = s_{c,i},$$

for all admissible transitions. This relation is unconditional and follows directly from the parity-controlled form of the Collatz rule.

The update of the octave parity $s_{a,i+1}$, however, is *not governed by a single universal identity*. Instead, it depends on the *affine form of the octave update* and on whether that update introduces a carry from higher-order-bits of the octave index.

More precisely, the octave index update can be written in normalized form as

$$A_{i+1} = (mA_i + c)/2^r.$$

Here $m \in \{1,3\}$ and c is the case-dependent offset determined by the base transition (Section 5). The exponent $r \geq 0$ is the number of factors of 2 extracted in the normalization step: After computing the raw affine numerator $mA_i + c$, one divides by $2^r$ so that the resulting base–octave representation (equations (20)–(21)) is restored. Equivalently, $r$ is the 2-adic valuation of the relevant numerator segment used to reconstitute $A_{i+1}$ in the canonical decomposition. The resulting parity $s_{a,i+1}$ satisfies one of two behaviors.

**Two admissible outcomes for $s_{a,i+1}$.** The finite parity descriptor Σ records only low-order binary information about $A_i$, namely $s_{a,i}$ (octave parity) and $s_{q,i}$ (the next binary digit of $A_i$ beyond $s_{a,i}$). When the octave index is updated by the affine rule and then renormalized into canonical base–octave form, higher-order bits of $A_i$ can introduce a carry across a power-of-two boundary. Consequently, for a fixed source descriptor $\Sigma_i$ and fixed base transition case, there are (in general) exactly two admissible outcomes for $s_{a,i+1}$ (and therefore for $s_{q,i+1}$). These two outcomes correspond to the carry/no-carry split enumerated in Appendix A1.

$$\text{(i) } s_{a,i+1} = s_{q,i},$$

$$\text{(ii) } s_{a,i+1} = 1 - s_{q,i}.$$

Operationally, for a given step the underlying integer update $h_i \to h_{i+1}$ fixes $A_{i+1}$ uniquely, hence it fixes $s_{a,i+1}$ uniquely. In the finite codebook representation (Appendix A1), the two admissible outcomes (i)–(ii) are listed explicitly as separate rows; an outcome ID label (or carry/no-carry label) is therefore only an index distinguishing these two rows. In computations, the correct branch is selected by matching the observed $s_{a,i+1}$ (computed from $h_{i+1}$) to the corresponding $s_{a,\,i+1}$ entry in Appendix A1.

These two outcomes are exactly the "carry/no-carry" split recorded in Appendix A1 in the $s_{a,i+1}$ and $s_{q,i+1}$ columns; for additional refinement see section 9.6. In contrast, the bits $s_{b,i}$ and $s_{c,i}$ track base-side parity structure (Section 7.1, equation (25)) and are governed by the unconditional base-parity identity; they are not affected by octave carry in $A$.

Thus, while the refined evolution is deterministic once the applicable affine/carry case is fixed, the octave-parity update cannot, in general, be expressed as a function of $\Sigma_i$ alone without reference to higher-order structure of $A_i$.

Appendix A1 enumerates all admissible refined transition cases, including the carry/branch outcomes that arise when higher-order bits of $A_i$ affect parity inheritance. In other words, the table should be read as a complete list of admissible case transitions $(A_i, \Sigma_i)\rightarrow(A_{i+1},\Sigma_{i+1})$, rather than as a single-valued map from $\Sigma_i$ alone.

**Remark (Determinism versus projection).** The Collatz update is a deterministic function on integers: given a current iterate $h_i$ (equivalently, given the pair $(B_i,A_i)$ in the base–octave decomposition), there is a unique successor $h_{i+1}$ (equivalently, a unique successor pair $(B_{i+1},A_{i+1})$). The refined parity descriptor $\Sigma_i$ is a finite projection of the full octave information carried by $A_i$. Because the affine octave update may introduce carries across binary boundaries, the next-step values of some displayed parity bits (e.g. $s_{a,i+1}$, $s_{q,i+1}$) need not be determined by the displayed tuple alone. This does not represent nondeterminism of the underlying iteration; rather, it reflects loss of information under projection. In Appendix A1 we therefore list admissible carry/branch outcomes explicitly as separate rows, so that each listed case corresponds to a unique successor state and affine octave update. This conditional structure suffices for the subsequent analysis: The convergence arguments depend only on the base transition routing and on monotone properties of the octave index, not on a universal parity-flip rule.

## 7.3 Illustration in the lowest octaves

To illustrate how the generalized Collatz rule and the base–octave decomposition operate concretely, Table 2 lists the values of $h_i$, together with the associated base index $B_i$, octave index $A_i$, octave parity par $(A_i)$, and resulting base value $B_{i+1}$, for integers $h = 1, \ldots, 16$ (the first two octaves).

This table serves a purely illustrative role: it demonstrates how familiar small integers are organized by the base–octave framework and how the selection rules act uniformly across both odd and even values. Higher octaves follow the same update rules and differ only in the scale parameter $A_i$.

### 7.4 Finite transition structure

The refined state description shows that the Collatz iteration induces a finite set of parity states coupled to an infinite family of octave levels. For each admissible refined state $\Sigma_i$, the successor $\Sigma_{i+1}$ is uniquely determined once the applicable affine/carry case is specified. For completeness, the full list of admissible parity-pathway transition cases—covering all combinations of base values, refined parity variables, and carry/branch outcomes—is given in Appendix A1. The appendix table functions as a reference codebook for these cases and is not required for understanding the derivations in Sections 5–7.

## 8. Growth and Decay Mechanisms Along Admissible Pathways

Sections 5–7 formulate the Collatz iteration as a parity-controlled dynamical system on the base–octave decomposition

$$h_i = B_i + 8(A_i - 1), B_i \in \{1, \dots, 8\},\ A_i \in \mathbb{N},$$

together with a refined parity state descriptor

$$\Sigma_i = (B_i, s_{b,i}, s_{c,i}, s_{a,i}, s_{q,i}).$$

In this representation, the evolution is deterministic: admissible states have unique successors, and the base-parity identity $s_{b,i+1} = s_{c,i}$ holds unconditionally, while the octave-parity update is conditional as described in Section 7.2. All admissible transitions of $\Sigma_i \to \Sigma_{i+1}$ are in Appendix A1. Representative trajectories grouped by base class are provided in Appendix B.

The purpose of the present section is to interpret the selection rule dynamics in terms of growth permitting versus decay-forcing channels. The goal here is mechanistic structure: to identify where large excursions can occur, where persistence is possible, and what structural constraints force eventual return to contractive regimes.

### 8.1 Base classes versus the integer terminal basin 1-2 cycle

A central distinction is the difference between a **base label** and an **integer value**. The base label $B$ = 1 does not mean $h$ =1 unless $A$ = 1. In general

$$B = 1 \Leftrightarrow h = 1 + 8(A - 1) = 8A - 7, B = 2 \Leftrightarrow h = 2 + 8(A - 1) = 8A - 6.$$

Thus, base classes $B = 1$ and $B = 2$ occur at arbitrarily high octaves and play an essential role for large iterates.

By contrast, the literal integer oscillation $1 \leftrightarrow 2$ occurs only when the trajectory reaches $A = 1$, i.e. when the iterate is $h = 1$ or $h = 2$ (integers 1 and 2), but base transitions B1 to B2, or B2 to B1 at higher octave levels ($h > 2$). This distinction is consistent with the "terminal component" viewpoint: termination corresponds to entry into the terminal strongly connected component containing the integer cycle $1 \leftrightarrow 2$, not to a single parity quadruple.

## 8.2 The {*B* = 1, *B* = 2} subsystem: mixed routing with strong octave contraction

The base classes $B = 1$ and $B = 2$ form a coupled subsystem that can occur at all octave levels. Although base-1 steps may increase scale locally, base-2 steps impose strong octave contraction (halving), and base-1 also routes into the junction base $B = 6$.

### Lemma 8.1 (Base-1 routing and octave update).

Let $h = 8A - 7$ (i.e. $B = 1$). Under the odd-update $h' = (3h + 1)/2$,

$$h' = 12A - 10.$$

Consequently,

$$B' = \begin{cases} 2, & A \text{ odd}, \\ 6, & A \text{ even}, \end{cases} \text{ and } A' = \begin{cases} 1 + \dfrac{3(A-1)}{2}, & A \text{ odd}, \\ \dfrac{3A-2}{2}, & A \text{ even}. \end{cases}$$

Thus a base-1 step may increase octave scale (approximately by a factor $3/2$), but it immediately routes either to base-2 or to the junction base-6.

### Lemma 8.2 (Base-2 halving mechanism).

Let $h = 8A - 6$ (i.e. $B = 2$). Under the even-update $h' = h/2$,

$$h' = 4A - 3.$$

Consequently,

$$B' = \begin{cases} 1, & A \text{ odd}, \\ 5, & A \text{ even}, \end{cases} \text{ and } A' = \begin{cases} \dfrac{A+1}{2}, & A \text{ odd}, \\ \dfrac{A}{2}, & A \text{ even}. \end{cases}$$

In particular, every visit to base $B = 2$ produces strong octave contraction (exact halving in even octaves and near-halving in odd octaves), even when $A$is large.

The selection rule - alternations $1 \leftrightarrow 2$ that occur at high octave levels represent base index - dynamics coupled to octave change; they are not an analogue of the terminal integer cycle. The {1,2} base subsystem is best viewed as a mixed routing mechanism in which base-2 steps enforce rapid octave contraction, while base-1 steps either feed base-2 or pass control to the junction base $6$.

### 8.3 Odd-base-gateways *B* = 3 and *B* = 5: local growth injection with rapid handoff to even bases

Odd-bases $B = 3$ and $B = 5$ do not support internal persistence (no $3 \to 3$ or $5 \to 5$ channels in the selection rules). Instead, they act as gateways: they may contribute local odd step-growth but then route deterministically into even bases.

Base $3$ routes to $5$ in odd octaves and to $1$ in even octaves; the subsequent pathways pass quickly into even bases via $5$ or via the {1,2] subsystem. Base $5$ routes to $8$ in odd octaves and to $4$ in even octaves, placing the trajectory immediately into the even base contraction network.

### Lemma 8.3a (Odd-base gateways are transient; no persistence outside $7 \to 7$).

Along any trajectory in base–octave form, every visit to base $B = 5$ is followed immediately by a visit to an even base, and every visit to base $B = 3$ is followed within at most two steps by a visit to an even base. In particular, bases 3 and 5 cannot support persistence or sustained growth: any growth they inject is forced to hand off rapidly into the even-base contraction network.

**Proof.** The selection rules (Section 5) give the complete list of admissible base transitions. For $B = 5$, the only possible outgoing transitions are $5 \to 8$ (odd octaves) and $5 \to 4$ (even octaves), both even bases.

For $B = 3$, the only possible outgoing transitions are $3 \to 5$ (odd octaves) and $3 \to 1$ (even octaves). If $3 \to 5$, then the very next step is $5 \to 4$ or $5 \to 8$, hence an even base is reached in one additional step. If $3 \to 1$, then from $B = 1$ the only outgoing transitions are $1 \to 2$ (odd octaves) or $1 \to 6$ (even octaves), both even bases, so an even base is reached in one additional step.

Thus $B = 5$ reaches an even base in one step and $B = 5$ reaches an even base in at most two steps. ■

This explains why excursions may include brief passages through 3 and 5 without enabling sustained growth: the permitted pathways force rapid handoff to contractive bases.

### 8.4 Base $B$ = 4: the principal even-base distributor (fork to 2 or to junction 6)

Base 4 is the central fork in the even-base network: In odd octaves, $4 \to 2$, frequently routing rapidly into the terminal basin on integer values (via 2 and 1); and in even octaves, $4 \to 6$, routing into the junction base 6.

Thus, trajectories reaching base 4 are quickly either driven toward the strongly contractive subsystem centered on base 2, or passed to base 6, where the system chooses between contraction via base 3 and potential delay via base 7.

### 8.5 Base $B$ = 6: the junction controlling access to delay

Base 6is the main junction controlling access to the delay channel:

In odd octaves, $6 \to 3$, which then routes through the gateway mechanism and returns to even bases; and in even octaves, $6 \to 7$, which may enter the base-7 persistence channel.

Accordingly, extended numerical excursions are structurally associated with pathways passing through the junction $6 \to 7$, because base-7 is the only base class supporting a persistence mechanism.

### 8.6 Base $B$ = 7: the delay channel and forced exit

Base 7 is special because it supports an internal even-octave-channel 7 → 7, whereas odd-octave visits eject immediately to $B$ = 3.

### Lemma 8.3 (Parity-gated base-7 channel).

Let $h = 8A - 1$ (i.e. $B = 7$). Under the odd-update $h' = (3h + 1)/2$,

$$h' = 12A - 1.$$

Consequently,

$$B' = \begin{cases} 7, & A \text{ even,} \\ 3, & A \text{ odd,} \end{cases} \quad \text{and if } A \text{ is even, } A' = \frac{3}{2}A.$$

Thus, the base 7 self-loop occurs only in even-octaves and scales the octave index by 3/2; once $A$ becomes odd, the next base transition is forced to exit 7 → 3.

### Proposition 8.4 (Bounded length of consecutive 7 → 7 runs).

Let $h_0 = 8A_0 - 1$. The number $t$ of consecutive base-7 self-loop steps 7 → 7 satisfies

$$t \leq \nu_2(A_0),$$

where $\nu_2(A_0)$ is the exponent of 2 dividing $A_0$. Indeed, each 7 → 7 step requires $A$ = even and updates $A \mapsto (3/2)A$, decreasing $\nu_2(A_0)$ by one. After at most $\nu_2(A_0)$ such steps, $A$ becomes odd and the next base transition is forced to exit via 7 → 3.

This quantifies the delay mechanism: base-7 persistence may be long when $A$ contains many factors of 2, but it cannot continue indefinitely.

### 8.6.1 Log-linear signature of growth episodes

### Lemma 8.6 (Exact log-linearity on consecutive odd-update runs).

Let $\{h_i\}_{i\geq 1}$ be a trajectory generated by equation (5). Suppose that for some index $n$ and some integer $t \geq 1$, the next $t$ updates are all of the odd type, i.e.

$$h_{k+1} = \frac{3h_k + 1}{2} \text{for all } k = n, n+1, \ldots, n+t-1.$$

Then the iterates satisfy the closed form

$$h_{n+t} + 1 = \left(\frac{3}{2}\right)^t (h_n + 1),$$

and therefore

$$\log_2(h_{n+t} + 1) = \log_2 (h_n + 1) + t \log_2\left(\frac{3}{2}\right)$$

In particular, on any interval of consecutive odd-updates, the graph of $\log_2(h_i+1)$ versus step index $i$ is exactly affine with slope $\log_2(3/2)$.

**Proof.** Define $x_k := h_k + 1$. Under the odd-update,

$$x_{k+1} = h_{k+1} + 1 = \frac{3h_k + 1}{2} + 1 = \frac{3h_k + 3}{2} = \frac{3}{2}(h_k + 1) = \frac{3}{2}x_k.$$

Iterating over $t$consecutive odd-updates gives $x_{n+t} = (3/2)^t x_n$, i.e. $h_{n+t} + 1 = (3/2)^t(h_n + 1)$. Taking $\log_2$ yields the stated affine relation. □

**Remark (Discrete "turning-point" detector in log coordinates).**

Let $y_i = \log_2(h_i + 1)$. On consecutive odd-updates one has $\Delta y_i = \log_2(3/2)$, hence $\Delta^2 y_i = 0$ throughout any growth segment where $\Delta y_i = y_{i+1} - y_i$ and $\Delta^2 y_i = \Delta y_{i+1} - \Delta y_i$, therefore, turning points between growth and decay correspond to indices where $\Delta^2 y_i \neq 0$, i.e., to discrete kinks in the piecewise-linear log profile. In the base 7 even-octave persistence channel, such kinks occur at the forced exit step when the entry level $A$ becomes odd; equivalently, after $\nu_2(A)$ growth steps the growth regime must terminate.

**Remark (Discrete log-average over growth intervals).**

For any interval of $t$ consecutive odd-updates, the identity

$$\log_2(h_{i+t} + 1) - \log_2(h_i + 1) = t \log_2(3/2)$$

is the exact telescoping sum of a constant discrete log-increment. Equivalently, the average per-step change of $\log_2(h + 1)$ over such an interval is exactly $\log_2(3/2)$. This is a finite summation identity (no limiting process and no continuous approximation are used).

Because the per-step log-growth on any growth interval is fixed while the admissible length of growth intervals is bounded by the channel structure, the cumulative logarithmic "height" gained in any single growth episode is explicitly bounded.

**Corollary 8.7 (Bounded log-gain per growth episode in the 7 → 7 channel).**

Consider a trajectory segment that remains in the base–octave channel $B = 7$ with even-octave persistence (the 7 → 7 channel). Let the segment begin at $h = 8A - 1$ with entry octave index $A = A_0$. Then the number $t$ of consecutive 7 → 7 steps is bounded by $t \leq \nu_2(A_0)$ (Proposition 8.4). Moreover, along this segment

$$\log_2(h + 1) = \log_2(8A) = 3 + \log_2(A)$$

increases exactly by $\log_2(3/2)$ per step (Lemma 8.6). Hence the total log-gain of any single 7 → 7 growth episode satisfies

$$\Delta \log_2(h + 1) \leq \nu_2(A_0) \log_2(3/2).$$

This provides an explicit limiter: long growth episodes require large 2-adic content in the entry octave index, and cannot persist indefinitely.

**Interpretation.** In log coordinates, growth intervals governed by consecutive odd-updates appear as straight-line segments in $\log_2(h_i + 1)$ with slope $\log_2(3/2)$ (Lemma 8.6), whereas decline intervals dominated by halving steps exhibit straight-line segments with slope $-1$ in $\log_2(h_i)$. The channel structure identifies where such growth segments can occur (Lemma 8.3) and bounds their duration (Proposition 8.4). When $\log_2(h_i + 1)$ is plotted against the iteration index, growth episodes appear as straight-line segments with slope $\log_2(3/2)$, while decline episodes governed by repeated halving appear as straight-line segments with slope $-1$; the duration of each growth segment is bounded by the channel structure (Lemma 8.3 and Proposition 8.4).

**Worked example (trajectory segment for $h_1$= 1639).**

Consider the generalized trajectory starting at $h_1 = 1639$ under the parity-controlled update (equation (5)). In this trajectory one encounters a base-7 persistence episode in an even octave: at some step *I* he iterate has the form $h_i = 8A_i - 1$ with $B_i = 7$ and $A_i$ even, so the base-7 channel 7 → 7 applies and the octave level updates by $A_{i+1} = \frac{3}{2}A_i$ while remaining in $B$= 7. For $h_1$ =1639, an explicit instance of this occurs when $h_i$ =1663, which corresponds to $B_i$ =7 and $A_i$ =208 (even). Along the subsequent base 7-persistence run, the iterate satisfies $h+1 = 8A$, hence each step multiplies $h+1$ by 3/2, so $\log_2(h+1)$ increases by exactly $\log_2(3/2)$ per step (Lemma 8.6). Moreover, because $\nu_2(208) = 4$, Proposition 8.4 predicts that exactly four consecutive 7 → 7 steps can occur before forced exit when $A$ becomes odd; indeed the run terminates precisely when the octave index loses its final factor of two, at which point the next base transition is 7 → 3 (Lemma 8.3 and Proposition 8.4).

In particular, the low-order binary structure of $A_i$ (e.g., $s_a, s_q$, and—if introduced as described later—$s_r$ as the next octave sub-parity bit) provides finite-state predictors for the available "growth budget", because each 7 → 7 step reduces $\nu_2(A)$ by one until forced exit.

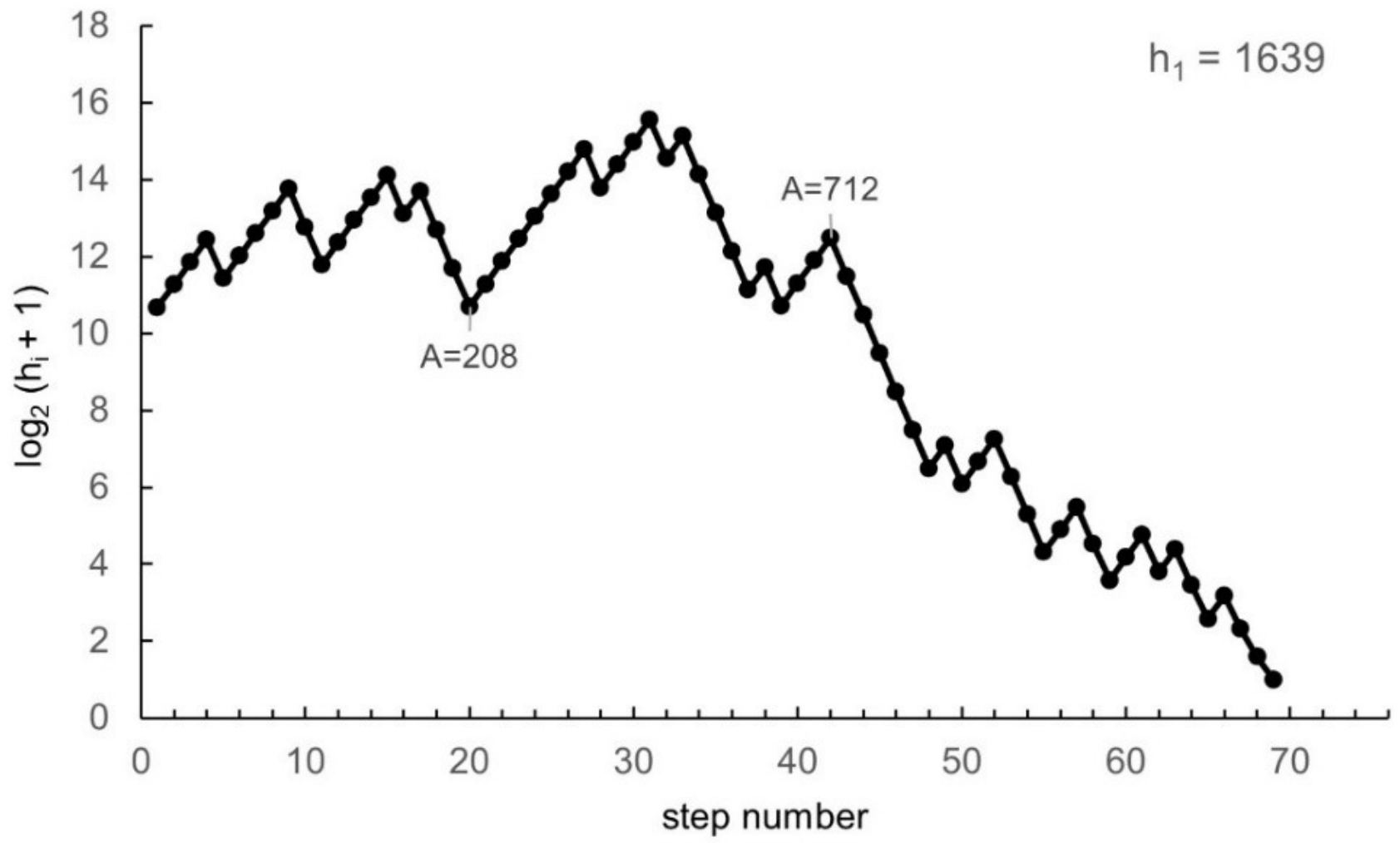

**Figure 2.** Piecewise linear growth and decay segments in log coordinates (illustrative example). Shown is the generalized Collatz trajectory starting at $h_i$ =1639 under the parity-controlled update (equation (5)), plotted as $y_i = \log_2(h_i+1)$ versus step index $i$. Consecutive odd-updates produce exact linear growth segments with slope $\log_2(3/2)$ in $y_i$ (Lemma 8.6), whereas halving steps produce linear

decay segments with slope -1. In this example, a growth episode is consistent with passage through the base-7 persistence channel 7→7 (Lemma 8.3), in which $A_{i+1} = 3A_i/2$ holds while $A_i$ remains even; the episode terminates at the forced exit step when $A_i$ becomes odd, consistent with the bound $t \leq \nu_2(A)$ (Proposition 8.4). The figure illustrates the log-linearity signatures and bounded duration of individual growth episodes implied by the channel structure. It is not evidence of monotone decline of peaks or of convergence: Trajectories may exhibit nonmonotone peak behavior (including secondary peaks). Section 9 completes the structural analysis using a return-map argument that establishes non-positive drift across persistence episodes.

## 8.7 Base $B = 8$: the collapse channel and forced drop to 4

Base 8 provides the strongest contraction mechanism. In even octaves it supports an internal channel $8 \to 8$; in odd octaves it drops to $B = 4$.

## Lemma 8.5 (Base-8 collapse channel).

Let $h = 8A$ (i.e. $B = 8$). Under the even-update $h' = h/2$,

$$h' = 4A.$$

Consequently,

$$B' = \begin{cases} 8, & A \text{ even,} \\ 4, & A \text{ odd,} \end{cases} \quad \text{and if } A \text{ is even, } A' = \frac{A}{2}.$$

Thus repeated 8 → 8 steps halve the octave index at each step and terminate once $A$ becomes odd, at which point the trajectory is forced to exit via 8 → 4 into the even-base contraction network. In particular, the length of an 8 → 8 contraction run in even octaves is *not* bounded a priori by the channel structure alone: it can be as long as the 2-adic valuation of the current octave index $A$ at entry, and hence can in principle be arbitrarily large across different visits/time scales.

Base-8 therefore acts as a collapse mechanism: it either transitions immediately to base-4 or remains temporarily in base 8 while rapidly decreasing octave level. This potential for long contraction segments complements the bounded-growth result for the base-7 channel: it explains how trajectories can exhibit extended descents in ν(A) after large excursions. The return-map analysis in section 9.6 shows that such regeneration cannot exceed the consumed budget, ruling out indefinite recurrence of secondary peaks.

### 8.8 Synthesis: where growth can occur and why it is constrained

The mechanism decomposition yields a coherent structural picture:

1. Persistent growth requires persistence in a growth-permitting regime**.** Among the base classes, the only internal persistence channel compatible with repeated expansive updates is the base 7 even-octave-channel 7 → 7.
2. That persistence is structurally bounded**.** Proposition 8.4 bounds the length of any consecutive 7 → 7 run by the 2-adic content of the octave index. Thus, long excursions are possible but must ultimately exit to the gateway bases and even-base-contraction network.
3. Odd bases 3 and 5 are transient gateways. They may contribute local growth under odd-updates but route quickly into even bases 4 and 8 (Lemma 8.3a).
4. Even bases dominate contraction. Base 2 enforces octave halving; base 8 enforces rapid halving through repeated 8 → 8 steps; base 4 distributes between direct contraction through base 2 and routing through junction 6.
5. The junction 6 controls access to delay. Entry into the bounded base-7 channel occurs primarily through 7 → 7 in even octaves; otherwise, 6 routes into the contraction gateway- pathway 6 → 3.

These structural constraints explain how large excursions can occur yet remain organized by a small number of admissible channel types. Appendix A1 provides the full parity pathway codebook for the deterministic state transitions, and Appendix B supplies representative numerical trajectories illustrating the resulting behavior. We now interpret the finite transition structure in terms of growth/decay channels.

The bounded-growth mechanism has a direct signature in log coordinates. On any run of consecutive odd-updates, $\log_2(h+1)$ is exactly affine in step number with slope $\log_2(3/2)$ (Lemma 8.6). Decline runs dominated by even-updates satisfy $\log_2(h_{i+1}) = \log_2(h_i) - 1$, i.e. slope $-1$. Thus, growth and decay intervals appear as straight-line segments in $\log_2(h+1)$ (growth) and $\log_2(h)$ (decay), with their durations constrained by the channel structure (notably Proposition 8.4).

## 8.9. Interaction with the net change constraint

Section 3 established the net-change identity (equation (13)), and Case 1 shows that any nontrivial finite cycle with $h_1 = h_m > 2$ would necessarily satisfy the corresponding zero–net-change constraint (equation (14)). Thus, exclusion of nontrivial cycles is a necessary requirement for global convergence. The net-change identity provides a global bookkeeping relation that any periodic orbit would have to satisfy; however, this identity alone does not exclude such orbits without additional structural input (e.g., constraints from the finite transition skeleton and channel structure developed later).

This eliminates finite mechanisms by which a sequence could avoid convergence (namely: nontrivial finite cycles and recurrence of integer values before entry into the terminal basin). Any non-terminating trajectory would therefore have to be infinite and strictly nonrepeating on integer values. When combined with the channel structure derived in Section 8, this observation sharply constrains possible escape behavior, but it does not by itself rule out all infinite trajectories.

Section 8 shows that persistent growth is structurally confined to the base-7 even-octave channel, and Proposition 8.4 proves that consecutive runs within this channel are finitely bounded by the 2-adic valuation of the octave index. Every exit from the growth channel routes the trajectory into base classes that enforce octave halving or rapid collapse. Because repeated halving necessarily contributes negative net change, indefinite avoidance of contraction would require repeated regeneration of growth potential without revisiting any previously attained value.

Consequently, any hypothetical nonterminating trajectory would have to realize an infinite, nonrepeating sequence of admissible states in which bounded growth episodes repeatedly compensate for cumulative halving losses, while also avoiding entry into the terminal basin. The return-map analysis in section 9.6 excludes this infinite regeneration scenario by proving that all return paths have non-positive 2-adic budget.

Thus, the convergence question is reduced to ruling out the existence of an infinite admissible regeneration path consistent with both the net-change constraint of Section 3 and the channel structure of Section 8.

## 9. State and Level Minimization

The preceding sections express the Collatz iteration in base–octave form

$$h_i = B_i + 8(A_i - 1), B_i \in \{1, \dots, 8\},\ A_i \in \mathbb{N},$$

and refine the dynamics using a finite parity-augmented state descriptor

$$\Sigma_i = (B_i, s_{b,i}, s_{c,i}, s_{a,i}, s_{q,i}).$$

As established in Sections 5–7, the update of $\Sigma_i$is deterministic once the full affine form of the update is specified. The complete set of admissible state transitions is finite and is enumerated exhaustively in Appendix A1. In particular, while certain parity inheritance relations hold conditionally (Section 7.2), the arguments below depend only on base routing and octave-level- monotonicity.

This section organizes the results of Sections 3–8 into a Lyapunov-style framework: a finite structural state space coupled to a logarithmic "height" coordinate on the octave level. Because trajectories can exhibit non-monotone behavior in this height coordinate (including secondary peaks), ν(*A*) is not a Lyapunov function per se. Instead, the analysis identifies a natural return-map quantity (the episode-to-episode return drift $\Delta\nu_{\text{return}}(k)$) whose non-positivity is established in section 9.6, suggesting eventual confinement in a contractive subsystem.

All statements in this section are made within the parity-controlled symbolic framework and admissibility constraints developed above.

**Guide to Appendix A1:** Appendix A is not a catalog of numerical examples but a transition codebook for the refined state description. Rows list admissible refined transition cases, including the carry/branch outcomes needed when higher-order bits of *A* affect parity inheritance. Given the current refined state information and the applicable case, one reads across the corresponding row to obtain the successor refined state and the associated affine octave update. In the convergence arguments of Section 9, the Appendix A1 is used to justify finite routing properties (e.g. which base transitions are possible and which are not), not to compute trajectories step-by-step.

## 9.1 Finite structural skeleton

Projecting away the octave index *A,* the Collatz dynamics induce a finite directed graph on the base classes B ∈ {1,...,8}, with parity-dependent edges. This graph captures the structural routing of trajectories independently of scale. For each admissible refined state Σ, the successor is uniquely determined once the applicable affine/carry case is specified; all admissible cases are enumerated in Appendix A1. Determinism here refers to uniqueness of the successor after fixing the relevant case; octave-parity inheritance may depend on higher-order structure of *A* and is handled in the codebook by listing the corresponding carry/branch outcomes.

A crucial point is that base labels do not correspond to small integers except at $A = 1$: transitions involving bases $B = 1$ and $B = 2$ occur at arbitrarily large magnitudes. Terminal integer behavior (the 1 ↔ 2 oscillation) arises only after entry into the terminal basin $A = 1$.

### One-step contraction predictor from the refined parity bits.

Although the base-transition-skeleton is formulated in terms of $(B_i, s_{a,i})$, the refined parity bit $s_{c,i}$ has a direct dynamical meaning: by the unconditional identity $s_{b,i+1} = s_{c,i}$, the value of $s_{c,i}$ determines the parity of the next base $B_{i+1}$. In particular, $s_{c,i} = 0$ implies that the next state lies in an even base class, hence in the contraction-dominant-part of the transition network. Combined with the octave parity gate $s_{a,i}$, this provides a deterministic finite "direction label" for whether the next step is forced into the even-base- contraction regime.

## 9.2 The octave index as a Lyapunov-style coordinate

### Even-base contraction

If $B_i$ is even, the update reduces to halving:

$$h_{i+1} = \frac{h_i}{2}.$$

In base–octave form this yields:

- for odd octaves: $A_{i+1} = 1 + (A_i - 1)/2$,
- for even octaves: $A_{i+1} = A_i/2$.

Thus, for all $A_i > 1$, every visit to an even base strictly decreases $A_i$, with equality only at $A_i = 1$ (terminal basin).

## Strong collapse channels

Base $B$ = 8 admits an internal halving channel 8 → 8 with $A_{i+1} = A_i/2$; once the octave parity flips, the trajectory exits into the even-base- contraction network. Base $B$ = 2 similarly enforces strong octave contraction. Hence, whenever the trajectory enters bases 2, 4, or 8, the octave index undergoes repeated halving and therefore decreases monotonically and efficiently.

## 9.2.1 Refinement: a base-2 Lyapunov coordinate for the octave index

To make the scale dynamics explicit, we introduce a logarithmic Lyapunov-style coordinate on the octave level,

$$\nu(A_i) := \log_2(A_i).$$

We define the stepwise logarithmic increment

$$\Delta\nu_i := \nu(A_i+1) - \nu(A_i) = \log_2(A_i+1/A_i).$$

This choice is natural in the present setting because the canonical remapping (equations (20)–(21)) repeatedly extracts powers of 2. In particular, whenever the octave update can be written in normalized form $A_{i+1}=(mA_i+c)/2^r$ (Section 7.2), the exponent r contributes exactly −r to $\Delta\nu_i$ (plus the contribution from the numerator). Thus, the increment $\Delta\nu i$ measures the net "bit gain" or "bit loss" in the octave index per operational step.

Over a finite trajectory segment of length N, we define the finite-time average log-drift in the octave coordinate,

$$\bar{\nu}N := (1/(N-1)) \sum i=1N-1\ \Delta\nu_i = (1/(N-1)) \log_2(AN/A_1).$$

The equality on the right is the exact telescoping identity; it shows that $\bar{\nu}$N can be computed either as a stepwise average or from endpoints. Conceptually, however, the stepwise form is useful because it separates trajectories into growth-permitting versus decay-forcing intervals. In particular, Section 8 identifies where positive-scale increments can occur (notably along the base-7 even-octave persistence channel) and proves that the length of such growth runs is bounded (Proposition 8.4), while Section 9.2 shows that visits to even bases enforce octave contraction. Thus the Lyapunov coordinate provides a compact way to track "ascent budget" and "descent enforcement" along admissible pathways.

Important caveat. A Lyapunov coordinate is not, by itself, a proof of convergence unless one establishes a *strict* decrease under all admissible nonterminal steps (or, more generally, a strict decrease under a suitable return-map). Here $\nu(A_i)$ does not decrease at every step because bounded growth episodes exist. The role of ν in this paper is therefore structural: it provides a quantitative coordinate that separates ascent (growth-permitting) and descent (decay-forcing) intervals and supports exact bookkeeping identities (telescoping sums). Section 9.5 identifies the natural return-map quantity $\Delta\nu_{return}(k)$ between successive base-7 persistence episodes (Lemma 9.8). Section 9.6 suggests that this return drift is non-positive; this is the key termination criterion in the return-map scheme.

## 9.3 The delaying mechanism: bounded base-7 persistence episodes (7→7 runs)

**Terminology note**. In what follows, a "base-7 persistence episode" (or "7→7 episode") means a maximal consecutive run of transitions with base index $B$ = 7 and even octave parity, i.e. repeated applications of the 7→7 selection rule; this is distinct from a "growth interval" in the sense of Lemma 8.6, which refers to an interval of consecutive odd-updates and may occur outside base 7.
The only admissible mechanism capable of sustaining growth of the octave index across consecutive steps is the base 7 even-octave persistence channel, as shown by the valuation analysis in Section 8.

If $h_i = 8A_i - 1$ (so $B_i = 7$) and $A_i$ is even, then the odd-update preserves base-7 and

$$A_{i+1} = \frac{3}{2}A_i.$$

However, each such step consumes one factor of 2 from $A_i$. Consequently, the number $t$ of consecutive 7→7 steps is bounded by

$$t \leq \nu_2(A_i).$$

This bound depends only on the entry level and not on the surrounding trajectory. Thus, the 7→7 mechanism is a finite delay channel**:** it can increase magnitude temporarily but cannot persist indefinitely for any fixed entry level.

## 9.4 Net change constraint and exclusion of recurrence

Section 3 establishes a global net-change identity expressing the final value $h_m$ as the initial value reduced by cumulative halving contributions and increased by cumulative odd contributions. Convergence to $h_m$ =1 is possible only if a strict net-reduction equal to $h_1$ −1 occurs.

A steady-state relation $h_m = h_1 > 2$ would correspond to a nontrivial closed cycle. The absence of such solutions implies that no integer value $h > 2$ can occur more than once along a trajectory without creating a forbidden recurrence.

As noted in section 8.9, the net-change identity alone does not exclude nontrivial cycles; the channel/routing argument below supplies the additional structural restrictions.

## 9.5 Finite regeneration of growth-permitting states

For clarity, throughout section 9.5 the term "growth episode" refers specifically to a base-7 persistence episode (a consecutive 7→7 run) as defined in section 9.3; it does not refer to an arbitrary run of consecutive odd-updates (Lemma 8.6).

### Lemma 9.5 (Return-to-7 forces an even-base contraction step).

Consider a trajectory written in base–octave form $h_i = B_i + 8(A_i - 1)$. Suppose the trajectory exits base 7 via the forced transition 7 → 3 (i.e. after the base-7 channel terminates). Then before the next occurrence of base 7, the trajectory must traverse at least one even base (hence at least one halving step), and an even base is reached within at most two steps after the exit 7→3. Moreover, any re-entry into base 7 must occur through base 6 via the even-octave transition 6 → 7.

Because $B_i$ even implies $h_i$ even (Lemma 5.1), every visit to an even base corresponds to an even iterate and hence triggers the halving update $h_{i+1}=h_i/2$ under the generalized rule (equation (5)); equivalently, even-base steps lie in the contraction regime summarized in Section 9.2.

**Proof.** The claim is a finite routing statement on the base transition graph, and it follows by exhaustive enumeration of admissible outgoing edges as listed in Appendix A1 (Table A1):

(i) **Forced exit from base 7.** Once a base-7 persistence episode ends (i.e. once the octave index becomes odd), the next admissible base transition is 7→3. There is no other outgoing transition from 7 in the nonpersistent regime.

(ii) **From base 3 only to bases 1 or 5.** The only outgoing transitions from base 3 listed in Table A1 are 3 → 1 and 3→5 (with the choice determined by octave parity).

(iii) **Even base reached within at most two steps after 7→3.** If 3→5 occurs, then the next base transition is 5→8, which is an even base. If 3→1 occurs, then the next base transition is either 1→2 (even base immediately) or 1→6. In the latter case, the subsequent transition from base 6 is either 6→3 or 6→7; in particular, any return toward base 7 proceeds through base 6. Thus, along every itinerary beginning with the forced exit 7→3, an even base is encountered within at most two further transitions (specifically via 5→8 or 1→2), and therefore between any exit from base 7 and any later return to base 7 there is at least one even-base step.

(iv) **Unique re-entry into base 7.** To enter base 7 from a different base, the destination must satisfy $B_{i+1}$=7. Inspecting Table A1, the only rows with destination base 7 are from base 6. Hence the only re-entry into base 7 is via 6→7, and it occurs in an even octave. ■

## Lemma 9.6 (Even-base steps strictly decrease octave index).

If $B_i$ is even and $A_i > 1$, then the update $h_{i+1} = h_i/2$ implies $A_{i+1} < A_i$ (Section 9.2, equations for odd/even octaves). ■

**Lemma 9.7 (Inter-episode contraction is unavoidable).**
Along any trajectory expressed in base–octave form, every base 7 persistence episode (maximal consecutive 7→7 run) is finite and terminates by forced exit 7→3. Moreover, between the end of one base 7 persistence episode and the beginning of the next base 7 persistence episode (if any), the trajectory must traverse at least one even-base state, hence at least one halving step on the underlying integer iterate.

**Proof.** The finiteness of each 7→7 episode is Proposition 8.4. The routing claim is the finite-graph statement of Lemma 9.5: After 7→3 the next bases are in {1,5}, and within at most two further transitions an even base is reached (via 5→8 or 1→2). ■

## Lemma 9.8 (Return-drift formulation in the log–octave coordinate).

Let $\nu(A_i) := \log_2(A_i)$ and $\Delta\nu_i := \nu(A_{i+1}) - \nu(A_i) = \log_2(A_{i+1}/A_i)$. Let $\tau(1)<\tau(2)<\ldots$ denote the start indices of successive base-7 persistence episodes (i.e. maximal consecutive 7→7 runs), as in 9.3.

**Notation.** The index k in $\tau(k)$ and $\Delta\nu_{return}(k)$ is an episode counter (k=1,2,...) and is unrelated to the decomposition indices $k_{b,i}$, $k_{a,i}$, $k_{q,i}$, $k_{c,i}$ introduced in Section 4.

$$\Delta\nu_{return}(k) := \sum_{i=\tau(k)}^{\tau(k+1)-1} \Delta\nu_i = \log_2(A_{\tau(k+1)}/A_{\tau(k)}).$$

Then $\Delta\nu_{return}(k)$ is exactly the net logarithmic scale change between successive base-7 persistence episodes. By Lemma 9.7, between $\tau(k)$ and $\tau(k+1)$ the trajectory must traverse at least one even base, hence incurs at least one strictly negative increment $\Delta\nu_i$ on that interval whenever $A > 1$.

However, the existence of a negative increment somewhere in the interval does *not* by itself imply $\Delta\nu_{return}(k)<0$; equivalently, it does not imply monotone decrease of the entrance level $A_{\tau(k)}$ across successive episodes.

**Remark.** Section 9.6 establishes that $\Delta\nu_{return}(k) \leq 0$ for all k, completing the Lyapunov/return-map style termination argument. The return-drift quantity is also the natural diagnostic in computational implementations because it is a telescoping sum of the stepwise increments $\Delta\nu_i$. and can be compared directly to the log-linear episode structure illustrated in Figure 2.

**Closing the return-map gap.** Sections 9.1–9.4 reduce global termination to controlling the return-map on the octave level between successive base 7 persistence episodes. Each persistence episode is finite (bounded by the 2-adic valuation of its entry octave index), and any re-entry into base 7 must traverse even-base states and therefore includes halving steps. The remaining task is to show that the inter-episode routing cannot regenerate the 2-adic "budget" consumed by persistence. We close this gap by refining the finite descriptor to a 128-state extension that includes one additional octave sub-parity bit and by exhaustively enumerating all admissible return paths from forced exit 7→3 to re-entry 6→7. For each admissible return path, we compute the net 2-adic budget and show it is non-positive; consequently, the episode-to-episode return drift in $\nu(A)=\log_2 A$ is non-positive. As a result, no admissible return path allows indefinite regeneration of persistence potential, and all admissible continuations eventually lie within the strictly contractive even-base network.

### 9.6.1. Return-map analysis and convergence structure

We formalize the return-map argument by refining the parity–base descriptor and controlling the evolution of the octave index between successive base-7 persistence episodes. We have (i) finite graph, (ii) reduce to simple paths, (iii) list all simple paths, (iv) compute budget, (v) conclude drift ≤ 0.

*Path Enumeration Procedure*: The 22 return paths catalogued in Table A2 were obtained by the following systematic enumeration on the 128-state transition graph G = (V, E).

*Graph construction*. Vertices V correspond to the 128 extended states $\Sigma' = (B, s_b, s_c, s_a, s_q, s_r)$; directed edges E encode the deterministic successor transitions listed in Table A1.

*Return subgraph*. We restrict attention to the subgraph G_R ⊂ G connecting the exit set $E_{73}$ = {states with $B$ = 7 and $s_a$ = 1, triggering the forced 7→3 transition} to the entry set $E_{67}$ = {states with $B$ = 6 and $s_a$ = 0, permitting the 6→7 re-entry}.

*Simple path enumeration*. A depth-first search (DFS) enumerates all simple paths (paths with no repeated vertices) from $E_{73}$ to $E_{67}$. The stopping rules are: (a) path reaches $E_{67}$ → record as valid return path; (b) path revisits a vertex → prune (cycle detected; by Remark 9.6.5a, cycles have non-positive weight); (c) no outgoing edges toward $E_{67}$ → prune (dead end).

*Budget computation*. For each enumerated path P, we compute Net_Budget(P) = $\Sigma$ $\nu_2$_generated(e) − $\Sigma$ $\nu_2$_consumed(e) − 1, where the −1 accounts for the mandatory entry cost at 6→7.

This procedure yields exactly 22 simple return paths. All satisfy Net_Budget(P) ≤ 0, as recorded in Table A2.

### Definition 9.6.1 (Extended parity state and branch index).

Extend the refined state $\Sigma = (B, s_{b,i}, s_{c,i}, s_{a,i}, s_{q,i})$ to

$\Sigma_i' := (B_i, s_{b,i}, s_{c,i}, s_{a,i}, s_{q,i}, s_{r,i})$

where $s_r$ = parity($k_q$), with $A$ written in dyadic form $A = 1 + 2k_a + 4k_q + 8k_r + \cdots$.

The bit $s_{r,i} \in \{0,1\}$ tracks the next level of 2-adic structure beyond sq. The extended states form a finite set of 128 admissible cases enumerated in Table A1. For a given integer iterate, the successor state is unique; the two admissible rows for each source state correspond to the two carry/no-carry outcomes determined by the actual value of $s_{a,i+1}$. The row index (OutcomeID) is a bookkeeping device and plays no role in the dynamics

### Definition 9.6.2 (2-adic budget)

For a transition e in the return subgraph, define:

$\nu_2$_consumed($e$) = 1    if $e$ is an even-base halving step ($B \in \{2,4,6,8\}$)

$\nu_2$_generated($e$) = 1    if $e$ regenerates a factor of 2 (odd step with $A$ odd)

For a return path $P$ from forced exit 7→3 to re-entry 6→7, define the net budget:

Net_Budget($P$) = $\Sigma\, \nu_2$_generated($e$) − $\Sigma\, \nu_2$_consumed($e$) − 1 ≤ 0

The −1 term accounts for the mandatory "entry cost": Re-entry 6→7 requires $s_a = 0$ (i.e., $A$ must be even at base 6), which consumes one factor of 2 from any regeneration.

### Lemma 9.6.4 (Bounded persistence).

Every base-7 persistence episode is finite, and its length is bounded above by the 2-adic valuation of the entry octave index $A_{entry}$.

**Proof.** Within base 7, each 7→7 step requires $s_a = 0$ ($A$ even) and produces $A_{i+1} = (3/2) \cdot A_i$, consuming one factor of 2 from the 2-adic decomposition. Because $\nu_2(A_{entry})$ is finite, only finitely many growth permitting steps are possible. ■

### Lemma 9.6.5 (Forced contraction on return).

Every return path from forced exit 7→3 to re-entry 6→7 contains at least one even-base halving step and hence includes at least one strictly contractive update of the octave index.

**Proof.** The finite transition skeleton (Table A1) shows that no admissible path connects 3 to 6 using only odd-base transitions. All return paths necessarily traverse even-base states where the update divides the octave index by two. ■

### Remark 9.6.5a (Completeness of simple paths).

The inter-episode dynamics define a finite directed graph on 128 states. Any cycle in the return subgraph has total weight ≤ 0: the 1↔2 and 8→8 loops are budget-neutral, whereas routes through 5→8 or 5→4 are strictly contractive. Therefore, the maximal net regeneration is attained on simple paths (no repeated vertices), and only 22 such paths exist from 7→3 to 6→7.

### Remark 9.6.5b (Exhaustiveness of the enumeration).

*Why 22 paths suffice*. The return subgraph is finite; it is therefore enough to evaluate simple paths, because any repeated state introduces a cycle whose net budget is non-positive.

1. *Finite state space*. The extended state space Σ' = $(B, s_b, s_c, s_a, s_q, s_r)$ contains exactly 128 states. The transition graph is fixed and finite. This is *not* sampling from infinitely many integers—it is enumerating paths in a finite directed graph.

2. *Simple paths suffice*: Any path with a repeated state contains a cycle. By Remark 9.6.5a, all cycles have weight ≤ 0. Therefore, adding a cycle to a simple path can only decrease Net_Budget. Conclusion: If all simple paths satisfy Net ≤ 0, then *all* paths (including infinitely long ones) satisfy Net ≤ 0.

3. *Complete enumeration*: The 22 paths in Table A2 are *all* simple paths from the exit set {7→3} to the re-entry set {6→7}. This is a standard finite graph enumeration problem, no different from proving a property holds for all edges in a graph.

4. *Trajectory aggregation at the state level*: Distinct integers $h_1 \neq h_2$ may correspond to the same extended state Σ' when they agree modulo sufficiently high powers of 2. In such cases, although the underlying integer values differ, the induced parity–base state is identical. From that point onward, state-level routing through the finite transition graph is the same, with differences persisting only in the octave index *A*. This aggregation highlights why the return-map analysis can be formulated entirely in terms of admissible state transitions and octave-level bookkeeping, without reference to the absolute size of the starting value. The convergence arguments depend only on properties of the finite state graph and the monotonic behavior of the octave coordinate along admissible paths.

A proof is therefore not by 'numerical evidence' but with a complete finite graph analysis. The enumeration *is* exhaustive in the precise mathematical sense.

### Proposition 9.6.6 (Non-positive return budget).

For every admissible return path *P* in the 128-state system,

Net_Budget(*P*) ≤ 0.

**Proof**. By exhaustive enumeration. Table A2 (Appendix) lists all 22 minimal simple return paths from 7→3 to 6→7. For each path, the table records $\nu_2$_consumed, $\nu_2$_max_gain, and Net_Budget (including the −1 entry cost).

The *entry cost* ensures that even paths achieving maximal regeneration satisfy Net_Budget ≤ 0. Paths containing cycles are strictly more contractive (Remark 9.6.5.a), and the enumeration is complete. ■

### Corollary 9.6.7 (No infinite persistence recurrence).

Base-7 persistence episodes cannot recur indefinitely along any trajectory.

**Proof.** By Lemma 9.6.4, each persistence episode consumes a positive amount of the $\nu_2$ budget. By Proposition 9.6.6, an inter-episode return path regenerates Net_Budget(*P*) ≤ 0. Thus,

$\nu_2$(A_next_entry) ≤ $\nu_2$(A_entry) − 1.

Hence, the available valuation strictly decreases across successive episodes, so infinite recurrence is impossible. ■

### Theorem 9.6.8 (Confinement to the Terminal Basin).

Every Collatz trajectory is eventually confined to the terminal basin {1,2}.

**Proof.** By Corollary 9.6.7, only finitely many growth-permitting base-7 persistence episodes can occur along any trajectory. After the final episode, the trajectory remains in the even-base subnetwork, where each update strictly decreases the octave index. Consequently, the trajectory reaches $A$ = 1 in finite time, at which point the induced dynamics enter the terminal 1↔2 basin, consistent with the Collatz conjecture. ■

## 9.7. Implementation check (non-essential).

The arguments in Sections 5–9 are entirely finite-state and combinatorial in nature. The convergence results follow from the structure of the parity-controlled transition graph and the exhaustive enumeration of admissible return paths in the refined 128-state system. No computational assumptions or numerical evidence are used in the structural arguments.

As an independent consistency check, a spreadsheet implementation of the refined state transitions and return-path enumeration was constructed. The implementation reproduces the admissible base routing and confirms that each enumerated return path satisfies the derived non-positive return-map bound. This computational check does not discover new cases, establish bounds, or replace analytical arguments; it merely audits the explicit finite enumeration carried out in Appendix A2 for internal consistency.

The spreadsheet plays no role in the derivation of the results and may be omitted without affecting any argument or conclusion. The spreadsheet implementation is available from the author upon request.

## 10. Conclusions

This paper reformulates the Collatz iteration as a parity-controlled dynamical system on a finite base–parity state space lifted across octave levels. In this representation, the dynamics decompose into a finite transition skeleton governing base evolution together with affine updates of an associated octave index. Refining the parity descriptor isolates the unique growth-permitting channel, corresponding to base-7 persistence at even-octave parity, and shows that such episodes are necessarily finite. Extending the state description to a finite 128-state system allows the admissible return paths between persistence episodes to be enumerated explicitly.

The critical observation is that re-entry into base-7 persistence requires an even octave index, imposing an entry constraint that excludes net regeneration of the 2-adic budget. The associated return-map analysis shows that the logarithmic octave coordinate exhibits non-positive drift across episode returns.

Therefore, base-7 persistence episodes cannot recur indefinitely. After the final such episode, all admissible continuations lie entirely within the contractive even-base subnetwork, and the subsequent evolution is confined to this regime. The resulting confinement suggests convergence to the terminal basin {1,2} within the parity-controlled framework developed here.

This approach emphasizes structural organization and finite-state mechanisms rather than ad hoc estimates, and demonstrates that parity-driven dynamical systems with bounded growth channels admit analysis via finite-state return-map methods, which may be applicable to related problems. This work illustrates how parity-controlled symbolic reductions can be leveraged to study integer recurrences with mixed expansive and contractive behavior.

***AI Assistance Disclosure**: Claude (Anthropic, 2024-2026) was used as a collaborative tool for: (1) structuring mathematical arguments supplied by the author, (2) drafting and editing manuscript text,* *and (3) creating spreadsheet implementations. The author retains full responsibility for all mathematical claims and proofs.*

***Funding Disclosure:*** *This research was conducted independently and did not receive any specific grant from funding agencies in the public, commercial, or not-for-profit sectors.*

## Appendix A1. Complete Parity Pathway Table

### A.1 Scope and Completeness

Appendix A functions as a complete transition codebook for the refined Collatz state introduced in Section 7. The refined parity pathway table consists of 128 admissible cases; Table A1 enumerates all of them. Each row corresponds to one admissible refined transition case and records the resulting successor refined state Σ', together with the associated base transition.

Completeness follows from combining the exhaustive base-transition derivation in Section 5 (Proposition 5.3) with the restricted form of the octave updates derived in Section 7.2. For every admissible transition, the octave index update can be written in affine form:

$$A(i+1) = (m \cdot A(i) + c) / 2^r$$

where $m \in \{1, 3\}$ and the offset c is determined by the base transition and octave parity (Sections 5 and 7.2). The exponent $r \geq 0$ records the number of powers of two removed in the normalization step. Carry/no-carry effects correspond to whether the relevant 2-adic valuation equals r or exceeds it by one, yielding two admissible outcomes recorded in the $s_a(i+1)$ and $s_q(i+1)$ columns of Table A1. Table A1 shows the essential transition data.

### A.2 Column Definitions for table A1.

Table A1 contains 22 columns, the column headers and their meanings are:

SOURCE STATE COLUMNS (columns 1-8):

1. StateID_128 Unique identifier encoding all state bits
2. B Base class (1-8), where $h = B + 8(A-1)$
3. $s_b$ Parity of $B$: $s_b = B \bmod 2$
4. $s_c$ Parity of $\lceil B/2 \rceil$: determines next base parity ($s_b(i+1) = s_c(i)$)
5. $s_a$ Parity of $A$: $s_a = A \bmod 2$ (0 = even, 1 = odd)
6. $s_q$ Next-level parity: $s_q = \lfloor A/2 \rfloor \bmod 2$
7. $s_r$ Third-level parity: $s_r = \lfloor A/4 \rfloor \bmod 2$
8. OutcomeID Row index (0 or 1) distinguishing two successor outcomes for the same source state. This is a *bookkeeping device* only; for any concrete integer, exactly one outcome is realized.

CLASSIFICATION COLUMNS (columns 9-11):

9. $v_2$_class 2-adic valuation class of A, determined by ($s_a$, $s_q$, $s_r$):

$s_a$=1 → $v_2$=0 ($A$ odd)
$s_a$=0, $s_q$=1 → $v_2$=1
$s_a$=0, $s_q$=0, $s_r$=1 → $v_2$=2
$s_a$=0, $s_q$=0, $s_r$=0 → $v_2 \geq 3$

10. max_persist   Maximum consecutive 7→7 steps possible from this state
11. drift_type   DECAY(-1), MIXED, or PERSIST(+0.585)

SUCCESSOR STATE COLUMNS (columns 12-16):
12. NextB   Successor base $B(i+1)$
13. next_$s_b$   Successor $s_b(i+1)$ [= $s_c(i)$, unconditionally]
14. next_$s_c$   Successor $s_c(i+1)$
15. next_$s_a$   Successor octave parity $s_a(i+1)$
16. next_$s_q$   Successor $s_q(i+1)$
17. next_$s_r$   Successor $s_r(i+1)$

FLAG COLUMNS (columns 18-20):
18. IsS7persist   TRUE if state is in B7 persistence set ($B$=7 AND $s_a$=0)
19. IsEntry67   TRUE if this is a 6→7 re-entry transition
20. IsExit73   TRUE if this is a forced 7→3 exit transition

BUDGET COLUMNS (columns 21-22):
21. $v_2$_consumed   1 if even-base halving step, else 0
22. $v_2$_possible_gain   1 if odd step can regenerate factor of 2, else 0

**A.3. How to read the parity pathway table.**

Rows are indexed by admissible refined cases, not by base value alone. Each row corresponds to one admissible refined transition case and specifies the corresponding successor refined state Σ'($i$+1). Where two rows share the same source tuple ($B$, $s_b$, $s_c$, $s_a$, $s_q$, $s_r$) but differ in the destination $s_a(i+1)$, $s_q(i+1)$, and/or $s_r(i+1)$, this reflects the two admissible carry/no-carry outcomes described in Section 7.2. The "OutcomeID" column (0 or 1) distinguishes these rows; it is a bookkeeping index only.

Important: For a given integer iterate, the successor state is *unique*. The two rows per source state represent the two *possible* outcomes; the *realized* outcome is determined by the actual value of $s_a(i+1)$ after computing the update. No admissible refined transition occurs outside the cases listed in Table A1; equivalently, Table A1 is exhaustive for the refined dynamics described in Sections 5–7.

All structural arguments in the main text rely solely on the base transitions and octave updates listed explicitly in the table.

## Appendix A1, Table A1. Complete Parity Pathway Table

| State ID (1) | $B_i$ (2) | $s_b$ (3) | $s_c$ (4) | $s_a$ (5) | $s_q$ (6) | $s_r$ (7) | OutCome ID (8) | $v_2$ Class (9) | Max Per. (10) | Drift Type (11) | next B (12) | next $s_b$ (13) | next $s_c$ (14) | next $s_a$ (15) | next $s_q$ (16) | next $s_r$ (17) | 7-7 (18) | Entry 6-7 (19) | Exit 7-3 (20) | v2 cons (21) | v2. gain (22) |
|---|---|---|---|---|---|---|---|---|---|---|---|---|---|---|---|---|---|---|---|---|---|
| B1_sb1_sc0_sa0_sq0_sr0_ID0 | 1 | 1 | 0 | 0 | 0 | 0 | 0 | $v_2 \geq 3$ | 3+ | MIXED | 6 | 0 | 1 | 0 | 0 | 0 | | | | 1 | 0 |
| B1_sb1_sc0_sa0_sq0_sr0_ID1 | 1 | 1 | 0 | 0 | 0 | 0 | 1 | $v_2 \geq 3$ | 3+ | MIXED | 6 | 0 | 1 | 0 | 0 | 1 | | | | 1 | 0 |
| B1_sb1_sc0_sa0_sq0_sr1_ID0 | 1 | 1 | 0 | 0 | 0 | 1 | 0 | $v_2=2$ | 2 | MIXED | 6 | 0 | 1 | 0 | 1 | 0 | | | | 1 | 0 |
| B1_sb1_sc0_sa0_sq0_sr1_ID1 | 1 | 1 | 0 | 0 | 0 | 1 | 1 | $v_2=2$ | 2 | MIXED | 6 | 0 | 1 | 0 | 1 | 1 | | | | 1 | 0 |
| B1_sb1_sc0_sa0_sq1_sr0_ID0 | 1 | 1 | 0 | 0 | 1 | 0 | 0 | $v_2=1$ | 1 | MIXED | 6 | 0 | 1 | 1 | 0 | 0 | | | | 1 | 0 |
| B1_sb1_sc0_sa0_sq1_sr0_ID1 | 1 | 1 | 0 | 0 | 1 | 0 | 1 | $v_2=1$ | 1 | MIXED | 6 | 0 | 1 | 1 | 0 | 1 | | | | 1 | 0 |
| B1_sb1_sc0_sa0_sq1_sr1_ID0 | 1 | 1 | 0 | 0 | 1 | 1 | 0 | $v_2=1$ | 1 | MIXED | 6 | 0 | 1 | 1 | 1 | 0 | | | | 1 | 0 |
| B1_sb1_sc0_sa0_sq1_sr1_ID1 | 1 | 1 | 0 | 0 | 1 | 1 | 1 | $v_2=1$ | 1 | MIXED | 6 | 0 | 1 | 1 | 1 | 1 | | | | 1 | 0 |
| B1_sb1_sc0_sa1_sq0_sr0_ID0 | 1 | 1 | 0 | 1 | 0 | 0 | 0 | $v_2=0$ | 0 | MIXED | 2 | 0 | 1 | 0 | 0 | 0 | | | | 0 | 1 |
| B1_sb1_sc0_sa1_sq0_sr0_ID1 | 1 | 1 | 0 | 1 | 0 | 0 | 1 | $v_2=0$ | 0 | MIXED | 2 | 0 | 1 | 1 | 1 | 1 | | | | 0 | 1 |
| B1_sb1_sc0_sa1_sq0_sr1_ID0 | 1 | 1 | 0 | 1 | 0 | 1 | 0 | $v_2=0$ | 0 | MIXED | 2 | 0 | 1 | 0 | 0 | 0 | | | | 0 | 1 |
| B1_sb1_sc0_sa1_sq0_sr1_ID1 | 1 | 1 | 0 | 1 | 0 | 1 | 1 | $v_2=0$ | 0 | MIXED | 2 | 0 | 1 | 1 | 1 | 1 | | | | 0 | 1 |
| B1_sb1_sc0_sa1_sq1_sr0_ID0 | 1 | 1 | 0 | 1 | 1 | 0 | 0 | $v_2=0$ | 0 | MIXED | 2 | 0 | 1 | 0 | 0 | 0 | | | | 0 | 1 |
| B1_sb1_sc0_sa1_sq1_sr0_ID1 | 1 | 1 | 0 | 1 | 1 | 0 | 1 | $v_2=0$ | 0 | MIXED | 2 | 0 | 1 | 1 | 1 | 1 | | | | 0 | 1 |
| B1_sb1_sc0_sa1_sq1_sr1_ID0 | 1 | 1 | 0 | 1 | 1 | 1 | 0 | $v_2=0$ | 0 | MIXED | 2 | 0 | 1 | 0 | 0 | 0 | | | | 0 | 1 |
| B1_sb1_sc0_sa1_sq1_sr1_ID1 | 1 | 1 | 0 | 1 | 1 | 1 | 1 | $v_2=0$ | 0 | MIXED | 2 | 0 | 1 | 1 | 1 | 1 | | | | 0 | 1 |
| B2_sb0_sc1_sa0_sq0_sr0_ID0 | 2 | 0 | 1 | 0 | 0 | 0 | 0 | $v_2 \geq 3$ | 3+ | DECAY(-1) | 5 | 1 | 0 | 0 | 0 | 0 | | | | 1 | 0 |
| B2_sb0_sc1_sa0_sq0_sr0_ID1 | 2 | 0 | 1 | 0 | 0 | 0 | 1 | $v_2 \geq 3$ | 3+ | DECAY(-1) | 5 | 1 | 0 | 0 | 0 | 1 | | | | 1 | 0 |
| B2_sb0_sc1_sa0_sq0_sr1_ID0 | 2 | 0 | 1 | 0 | 0 | 1 | 0 | $v_2=2$ | 2 | DECAY(-1) | 5 | 1 | 0 | 0 | 1 | 0 | | | | 1 | 0 |
| B2_sb0_sc1_sa0_sq0_sr1_ID1 | 2 | 0 | 1 | 0 | 0 | 1 | 1 | $v_2=2$ | 2 | DECAY(-1) | 5 | 1 | 0 | 0 | 1 | 1 | | | | 1 | 0 |
| B2_sb0_sc1_sa0_sq1_sr0_ID0 | 2 | 0 | 1 | 0 | 1 | 0 | 0 | $v_2=1$ | 1 | DECAY(-1) | 5 | 1 | 0 | 1 | 0 | 0 | | | | 1 | 0 |
| B2_sb0_sc1_sa0_sq1_sr0_ID1 | 2 | 0 | 1 | 0 | 1 | 0 | 1 | $v_2=1$ | 1 | DECAY(-1) | 5 | 1 | 0 | 1 | 0 | 1 | | | | 1 | 0 |
| B2_sb0_sc1_sa0_sq1_sr1_ID0 | 2 | 0 | 1 | 0 | 1 | 1 | 0 | $v_2=1$ | 1 | DECAY(-1) | 5 | 1 | 0 | 1 | 1 | 0 | | | | 1 | 0 |
| B2_sb0_sc1_sa0_sq1_sr1_ID1 | 2 | 0 | 1 | 0 | 1 | 1 | 1 | $v_2=1$ | 1 | DECAY(-1) | 5 | 1 | 0 | 1 | 1 | 1 | | | | 1 | 0 |
| B2_sb0_sc1_sa1_sq0_sr0_ID0 | 2 | 0 | 1 | 1 | 0 | 0 | 0 | $v_2=0$ | 0 | DECAY(-1) | 1 | 1 | 0 | 0 | 0 | 0 | | | | 0 | 0 |
| B2_sb0_sc1_sa1_sq0_sr0_ID1 | 2 | 0 | 1 | 1 | 0 | 0 | 1 | $v_2=0$ | 0 | DECAY(-1) | 1 | 1 | 0 | 1 | 1 | 1 | | | | 0 | 1 |
| B2_sb0_sc1_sa1_sq0_sr1_ID0 | 2 | 0 | 1 | 1 | 0 | 1 | 0 | $v_2=0$ | 0 | DECAY(-1) | 1 | 1 | 0 | 0 | 0 | 0 | | | | 0 | 0 |
| B2_sb0_sc1_sa1_sq0_sr1_ID1 | 2 | 0 | 1 | 1 | 0 | 1 | 1 | $v_2=0$ | 0 | DECAY(-1) | 1 | 1 | 0 | 1 | 1 | 1 | | | | 0 | 1 |
| B2_sb0_sc1_sa1_sq1_sr0_ID0 | 2 | 0 | 1 | 1 | 1 | 0 | 0 | $v_2=0$ | 0 | DECAY(-1) | 1 | 1 | 0 | 0 | 0 | 0 | | | | 0 | 0 |
| B2_sb0_sc1_sa1_sq1_sr0_ID1 | 2 | 0 | 1 | 1 | 1 | 0 | 1 | $v_2=0$ | 0 | DECAY(-1) | 1 | 1 | 0 | 1 | 1 | 1 | | | | 0 | 1 |
| B2_sb0_sc1_sa1_sq1_sr1_ID0 | 2 | 0 | 1 | 1 | 1 | 1 | 0 | $v_2=0$ | 0 | DECAY(-1) | 1 | 1 | 0 | 0 | 0 | 0 | | | | 0 | 0 |
| B2_sb0_sc1_sa1_sq1_sr1_ID1 | 2 | 0 | 1 | 1 | 1 | 1 | 1 | $v_2=0$ | 0 | DECAY(-1) | 1 | 1 | 0 | 1 | 1 | 1 | | | | 0 | 1 |
| B3_sb1_sc1_sa0_sq0_sr0_ID0 | 3 | 1 | 1 | 0 | 0 | 0 | 0 | $v_2 \geq 3$ | 3+ | MIXED | 1 | 1 | 0 | 0 | 0 | 0 | | | | 1 | 0 |
| B3_sb1_sc1_sa0_sq0_sr0_ID1 | 3 | 1 | 1 | 0 | 0 | 0 | 1 | $v_2 \geq 3$ | 3+ | MIXED | 1 | 1 | 0 | 0 | 0 | 1 | | | | 1 | 0 |
| B3_sb1_sc1_sa0_sq0_sr1_ID0 | 3 | 1 | 1 | 0 | 0 | 1 | 0 | $v_2=2$ | 2 | MIXED | 1 | 1 | 0 | 0 | 1 | 0 | | | | 1 | 0 |
| B3_sb1_sc1_sa0_sq0_sr1_ID1 | 3 | 1 | 1 | 0 | 0 | 1 | 1 | $v_2=2$ | 2 | MIXED | 1 | 1 | 0 | 0 | 1 | 1 | | | | 1 | 0 |
| B3_sb1_sc1_sa0_sq1_sr0_ID0 | 3 | 1 | 1 | 0 | 1 | 0 | 0 | $v_2=1$ | 1 | MIXED | 1 | 1 | 0 | 1 | 0 | 0 | | | | 1 | 0 |

| State ID (1) | $B_i$ (2) | $s_b$ (3) | $s_c$ (4) | $s_a$ (5) | $s_q$ (6) | $s_r$ (7) | OutCome ID (8) | $v_2$ Class (9) | Max Per. (10) | Drift Type (11) | next B (12) | next $s_b$ (13) | next $s_c$ (14) | next $s_a$ (15) | next $s_q$ (16) | next $s_r$ (17) | 7-7 (18) | Entry 6-7 (19) | Exit 7-3 (20) | v2 cons (21) | v2. gain (22) |
|---|---|---|---|---|---|---|---|---|---|---|---|---|---|---|---|---|---|---|---|---|---|
| B3_sb1_sc1_sa0_sq1_sr0_ID1 | 3 | 1 | 1 | 0 | 1 | 0 | 1 | $v_2$=1 | 1 | MIXED | 1 | 1 | 0 | 1 | 0 | 1 | | | | 1 | 0 |
| B3_sb1_sc1_sa0_sq1_sr1_ID0 | 3 | 1 | 1 | 0 | 1 | 1 | 0 | $v_2$=1 | 1 | MIXED | 1 | 1 | 0 | 1 | 1 | 0 | | | | 1 | 0 |
| B3_sb1_sc1_sa0_sq1_sr1_ID1 | 3 | 1 | 1 | 0 | 1 | 1 | 1 | $v_2$=1 | 1 | MIXED | 1 | 1 | 0 | 1 | 1 | 1 | | | | 1 | 0 |
| B3_sb1_sc1_sa1_sq0_sr0_ID0 | 3 | 1 | 1 | 1 | 0 | 0 | 0 | $v_2$=0 | 0 | MIXED | 5 | 1 | 0 | 0 | 0 | 0 | | | | 0 | 1 |
| B3_sb1_sc1_sa1_sq0_sr0_ID1 | 3 | 1 | 1 | 1 | 0 | 0 | 1 | $v_2$=0 | 0 | MIXED | 5 | 1 | 0 | 1 | 1 | 1 | | | | 0 | 1 |
| B3_sb1_sc1_sa1_sq0_sr1_ID0 | 3 | 1 | 1 | 1 | 0 | 1 | 0 | $v_2$=0 | 0 | MIXED | 5 | 1 | 0 | 0 | 0 | 0 | | | | 0 | 1 |
| B3_sb1_sc1_sa1_sq0_sr1_ID1 | 3 | 1 | 1 | 1 | 0 | 1 | 1 | $v_2$=0 | 0 | MIXED | 5 | 1 | 0 | 1 | 1 | 1 | | | | 0 | 1 |
| B3_sb1_sc1_sa1_sq1_sr0_ID0 | 3 | 1 | 1 | 1 | 1 | 0 | 0 | $v_2$=0 | 0 | MIXED | 5 | 1 | 0 | 0 | 0 | 0 | | | | 0 | 1 |
| B3_sb1_sc1_sa1_sq1_sr0_ID1 | 3 | 1 | 1 | 1 | 1 | 0 | 1 | $v_2$=0 | 0 | MIXED | 5 | 1 | 0 | 1 | 1 | 1 | | | | 0 | 1 |
| B3_sb1_sc1_sa1_sq1_sr1_ID0 | 3 | 1 | 1 | 1 | 1 | 1 | 0 | $v_2$=0 | 0 | MIXED | 5 | 1 | 0 | 0 | 0 | 0 | | | | 0 | 1 |
| B3_sb1_sc1_sa1_sq1_sr1_ID1 | 3 | 1 | 1 | 1 | 1 | 1 | 1 | $v_2$=0 | 0 | MIXED | 5 | 1 | 0 | 1 | 1 | 1 | | | | 0 | 1 |
| B4_sb0_sc0_sa0_sq0_sr0_ID0 | 4 | 0 | 0 | 0 | 0 | 0 | 0 | $v_2$≥3 | 3+ | DECAY(-1) | 6 | 0 | 1 | 0 | 0 | 0 | | | | 1 | 0 |
| B4_sb0_sc0_sa0_sq0_sr0_ID1 | 4 | 0 | 0 | 0 | 0 | 0 | 1 | $v_2$≥3 | 3+ | DECAY(-1) | 6 | 0 | 1 | 0 | 0 | 1 | | | | 1 | 0 |
| B4_sb0_sc0_sa0_sq0_sr1_ID0 | 4 | 0 | 0 | 0 | 0 | 1 | 0 | $v_2$=2 | 2 | DECAY(-1) | 6 | 0 | 1 | 0 | 1 | 0 | | | | 1 | 0 |
| B4_sb0_sc0_sa0_sq0_sr1_ID1 | 4 | 0 | 0 | 0 | 0 | 1 | 1 | $v_2$=2 | 2 | DECAY(-1) | 6 | 0 | 1 | 0 | 1 | 1 | | | | 1 | 0 |
| B4_sb0_sc0_sa0_sq1_sr0_ID0 | 4 | 0 | 0 | 0 | 1 | 0 | 0 | $v_2$=1 | 1 | DECAY(-1) | 6 | 0 | 1 | 1 | 0 | 0 | | | | 1 | 0 |
| B4_sb0_sc0_sa0_sq1_sr0_ID1 | 4 | 0 | 0 | 0 | 1 | 0 | 1 | $v_2$=1 | 1 | DECAY(-1) | 6 | 0 | 1 | 1 | 0 | 1 | | | | 1 | 0 |
| B4_sb0_sc0_sa0_sq1_sr1_ID0 | 4 | 0 | 0 | 0 | 1 | 1 | 0 | $v_2$=1 | 1 | DECAY(-1) | 6 | 0 | 1 | 1 | 1 | 0 | | | | 1 | 0 |
| B4_sb0_sc0_sa0_sq1_sr1_ID1 | 4 | 0 | 0 | 0 | 1 | 1 | 1 | $v_2$=1 | 1 | DECAY(-1) | 6 | 0 | 1 | 1 | 1 | 1 | | | | 1 | 0 |
| B4_sb0_sc0_sa1_sq0_sr0_ID0 | 4 | 0 | 0 | 1 | 0 | 0 | 0 | $v_2$=0 | 0 | DECAY(-1) | 2 | 0 | 1 | 0 | 0 | 0 | | | | 0 | 0 |
| B4_sb0_sc0_sa1_sq0_sr0_ID1 | 4 | 0 | 0 | 1 | 0 | 0 | 1 | $v_2$=0 | 0 | DECAY(-1) | 2 | 0 | 1 | 1 | 1 | 1 | | | | 0 | 1 |
| B4_sb0_sc0_sa1_sq0_sr1_ID0 | 4 | 0 | 0 | 1 | 0 | 1 | 0 | $v_2$=0 | 0 | DECAY(-1) | 2 | 0 | 1 | 0 | 0 | 0 | | | | 0 | 0 |
| B4_sb0_sc0_sa1_sq0_sr1_ID1 | 4 | 0 | 0 | 1 | 0 | 1 | 1 | $v_2$=0 | 0 | DECAY(-1) | 2 | 0 | 1 | 1 | 1 | 1 | | | | 0 | 1 |
| B4_sb0_sc0_sa1_sq1_sr0_ID0 | 4 | 0 | 0 | 1 | 1 | 0 | 0 | $v_2$=0 | 0 | DECAY(-1) | 2 | 0 | 1 | 0 | 0 | 0 | | | | 0 | 0 |
| B4_sb0_sc0_sa1_sq1_sr0_ID1 | 4 | 0 | 0 | 1 | 1 | 0 | 1 | $v_2$=0 | 0 | DECAY(-1) | 2 | 0 | 1 | 1 | 1 | 1 | | | | 0 | 1 |
| B4_sb0_sc0_sa1_sq1_sr1_ID0 | 4 | 0 | 0 | 1 | 1 | 1 | 0 | $v_2$=0 | 0 | DECAY(-1) | 2 | 0 | 1 | 0 | 0 | 0 | | | | 0 | 0 |
| B4_sb0_sc0_sa1_sq1_sr1_ID1 | 4 | 0 | 0 | 1 | 1 | 1 | 1 | $v_2$=0 | 0 | DECAY(-1) | 2 | 0 | 1 | 1 | 1 | 1 | | | | 0 | 1 |
| B5_sb1_sc0_sa0_sq0_sr0_ID0 | 5 | 1 | 0 | 0 | 0 | 0 | 0 | $v_2$≥3 | 3+ | MIXED | 4 | 0 | 0 | 0 | 0 | 0 | | | | 1 | 0 |
| B5_sb1_sc0_sa0_sq0_sr0_ID1 | 5 | 1 | 0 | 0 | 0 | 0 | 1 | $v_2$≥3 | 3+ | MIXED | 4 | 0 | 0 | 0 | 0 | 1 | | | | 1 | 0 |
| B5_sb1_sc0_sa0_sq0_sr1_ID0 | 5 | 1 | 0 | 0 | 0 | 1 | 0 | $v_2$=2 | 2 | MIXED | 4 | 0 | 0 | 0 | 1 | 0 | | | | 1 | 0 |
| B5_sb1_sc0_sa0_sq0_sr1_ID1 | 5 | 1 | 0 | 0 | 0 | 1 | 1 | $v_2$=2 | 2 | MIXED | 4 | 0 | 0 | 0 | 1 | 1 | | | | 1 | 0 |
| B5_sb1_sc0_sa0_sq1_sr0_ID0 | 5 | 1 | 0 | 0 | 1 | 0 | 0 | $v_2$=1 | 1 | MIXED | 4 | 0 | 0 | 1 | 0 | 0 | | | | 1 | 0 |
| B5_sb1_sc0_sa0_sq1_sr0_ID1 | 5 | 1 | 0 | 0 | 1 | 0 | 1 | $v_2$=1 | 1 | MIXED | 4 | 0 | 0 | 1 | 0 | 1 | | | | 1 | 0 |
| B5_sb1_sc0_sa0_sq1_sr1_ID0 | 5 | 1 | 0 | 0 | 1 | 1 | 0 | $v_2$=1 | 1 | MIXED | 4 | 0 | 0 | 1 | 1 | 0 | | | | 1 | 0 |
| B5_sb1_sc0_sa0_sq1_sr1_ID1 | 5 | 1 | 0 | 0 | 1 | 1 | 1 | $v_2$=1 | 1 | MIXED | 4 | 0 | 0 | 1 | 1 | 1 | | | | 1 | 0 |
| B5_sb1_sc0_sa1_sq0_sr0_ID0 | 5 | 1 | 0 | 1 | 0 | 0 | 0 | $v_2$=0 | 0 | MIXED | 8 | 0 | 0 | 0 | 0 | 0 | | | | 0 | 1 |
| B5_sb1_sc0_sa1_sq0_sr0_ID1 | 5 | 1 | 0 | 1 | 0 | 0 | 1 | $v_2$=0 | 0 | MIXED | 8 | 0 | 0 | 1 | 1 | 1 | | | | 0 | 1 |
| B5_sb1_sc0_sa1_sq0_sr1_ID0 | 5 | 1 | 0 | 1 | 0 | 1 | 0 | $v_2$=0 | 0 | MIXED | 8 | 0 | 0 | 0 | 0 | 0 | | | | 0 | 1 |
| B5_sb1_sc0_sa1_sq0_sr1_ID1 | 5 | 1 | 0 | 1 | 0 | 1 | 1 | $v_2$=0 | 0 | MIXED | 8 | 0 | 0 | 1 | 1 | 1 | | | | 0 | 1 |
| B5_sb1_sc0_sa1_sq1_sr0_ID0 | 5 | 1 | 0 | 1 | 1 | 0 | 0 | $v_2$=0 | 0 | MIXED | 8 | 0 | 0 | 0 | 0 | 0 | | | | 0 | 1 |

| State ID (1) | $B_i$ (2) | $s_b$ (3) | $s_c$ (4) | $s_a$ (5) | $s_q$ (6) | $s_r$ (7) | OutCome ID (8) | $v_2$ Class (9) | Max Per. (10) | Drift Type (11) | next B (12) | next $s_b$ (13) | next $s_c$ (14) | next $s_a$ (15) | next $s_q$ (16) | next $s_r$ (17) | 7-7 (18) | Entry 6-7 (19) | Exit 7-3 (20) | v2 cons (21) | v2. gain (22) |
|---|---|---|---|---|---|---|---|---|---|---|---|---|---|---|---|---|---|---|---|---|---|
| B5_sb1_sc0_sa1_sq1_sr0_ID1 | 5 | 1 | 0 | 1 | 1 | 0 | 1 | $v_2=0$ | 0 | MIXED | 8 | 0 | 0 | 1 | 1 | 1 | | | | 0 | 1 |
| B5_sb1_sc0_sa1_sq1_sr1_ID0 | 5 | 1 | 0 | 1 | 1 | 1 | 0 | $v_2=0$ | 0 | MIXED | 8 | 0 | 0 | 0 | 0 | 0 | | | | 0 | 1 |
| B5_sb1_sc0_sa1_sq1_sr1_ID1 | 5 | 1 | 0 | 1 | 1 | 1 | 1 | $v_2=0$ | 0 | MIXED | 8 | 0 | 0 | 1 | 1 | 1 | | | | 0 | 1 |
| B6_sb0_sc1_sa0_sq0_sr0_ID0 | 6 | 0 | 1 | 0 | 0 | 0 | 0 | $v_2\geq3$ | 3+ | DECAY(-1) | 7 | 1 | 1 | 0 | 0 | 0 | | ENTRY | | 1 | 0 |
| B6_sb0_sc1_sa0_sq0_sr0_ID1 | 6 | 0 | 1 | 0 | 0 | 0 | 1 | $v_2\geq3$ | 3+ | DECAY(-1) | 7 | 1 | 1 | 0 | 0 | 1 | | ENTRY | | 1 | 0 |
| B6_sb0_sc1_sa0_sq0_sr1_ID0 | 6 | 0 | 1 | 0 | 0 | 1 | 0 | $v_2=2$ | 2 | DECAY(-1) | 7 | 1 | 1 | 0 | 1 | 0 | | ENTRY | | 1 | 0 |
| B6_sb0_sc1_sa0_sq0_sr1_ID1 | 6 | 0 | 1 | 0 | 0 | 1 | 1 | $v_2=2$ | 2 | DECAY(-1) | 7 | 1 | 1 | 0 | 1 | 1 | | ENTRY | | 1 | 0 |
| B6_sb0_sc1_sa0_sq1_sr0_ID0 | 6 | 0 | 1 | 0 | 1 | 0 | 0 | $v_2=1$ | 1 | DECAY(-1) | 7 | 1 | 1 | 1 | 0 | 0 | | ENTRY | | 1 | 0 |
| B6_sb0_sc1_sa0_sq1_sr0_ID1 | 6 | 0 | 1 | 0 | 1 | 0 | 1 | $v_2=1$ | 1 | DECAY(-1) | 7 | 1 | 1 | 1 | 0 | 1 | | ENTRY | | 1 | 0 |
| B6_sb0_sc1_sa0_sq1_sr1_ID0 | 6 | 0 | 1 | 0 | 1 | 1 | 0 | $v_2=1$ | 1 | DECAY(-1) | 7 | 1 | 1 | 1 | 1 | 0 | | ENTRY | | 1 | 0 |
| B6_sb0_sc1_sa0_sq1_sr1_ID1 | 6 | 0 | 1 | 0 | 1 | 1 | 1 | $v_2=1$ | 1 | DECAY(-1) | 7 | 1 | 1 | 1 | 1 | 1 | | ENTRY | | 1 | 0 |
| B6_sb0_sc1_sa1_sq0_sr0_ID0 | 6 | 0 | 1 | 1 | 0 | 0 | 0 | $v_2=0$ | 0 | DECAY(-1) | 3 | 1 | 1 | 0 | 0 | 0 | | | | 0 | 0 |
| B6_sb0_sc1_sa1_sq0_sr0_ID1 | 6 | 0 | 1 | 1 | 0 | 0 | 1 | $v_2=0$ | 0 | DECAY(-1) | 3 | 1 | 1 | 1 | 1 | 1 | | | | 0 | 1 |
| B6_sb0_sc1_sa1_sq0_sr1_ID0 | 6 | 0 | 1 | 1 | 0 | 1 | 0 | $v_2=0$ | 0 | DECAY(-1) | 3 | 1 | 1 | 0 | 0 | 0 | | | | 0 | 0 |
| B6_sb0_sc1_sa1_sq0_sr1_ID1 | 6 | 0 | 1 | 1 | 0 | 1 | 1 | $v_2=0$ | 0 | DECAY(-1) | 3 | 1 | 1 | 1 | 1 | 1 | | | | 0 | 1 |
| B6_sb0_sc1_sa1_sq1_sr0_ID0 | 6 | 0 | 1 | 1 | 1 | 0 | 0 | $v_2=0$ | 0 | DECAY(-1) | 3 | 1 | 1 | 0 | 0 | 0 | | | | 0 | 0 |
| B6_sb0_sc1_sa1_sq1_sr0_ID1 | 6 | 0 | 1 | 1 | 1 | 0 | 1 | $v_2=0$ | 0 | DECAY(-1) | 3 | 1 | 1 | 1 | 1 | 1 | | | | 0 | 1 |
| B6_sb0_sc1_sa1_sq1_sr1_ID0 | 6 | 0 | 1 | 1 | 1 | 1 | 0 | $v_2=0$ | 0 | DECAY(-1) | 3 | 1 | 1 | 0 | 0 | 0 | | | | 0 | 0 |
| B6_sb0_sc1_sa1_sq1_sr1_ID1 | 6 | 0 | 1 | 1 | 1 | 1 | 1 | $v_2=0$ | 0 | DECAY(-1) | 3 | 1 | 1 | 1 | 1 | 1 | | | | 0 | 1 |
| B7_sb1_sc1_sa0_sq0_sr0_ID0 | 7 | 1 | 1 | 0 | 0 | 0 | 0 | $v_2\geq3$ | 3+ | GROW(+0.585) | 7 | 1 | 1 | 0 | 0 | 0 | YES | | | 1 | 0 |
| B7_sb1_sc1_sa0_sq0_sr0_ID1 | 7 | 1 | 1 | 0 | 0 | 0 | 1 | $v_2\geq3$ | 3+ | GROW(+0.585) | 7 | 1 | 1 | 0 | 0 | 1 | YES | | | 1 | 0 |
| B7_sb1_sc1_sa0_sq0_sr1_ID0 | 7 | 1 | 1 | 0 | 0 | 1 | 0 | $v_2=2$ | 2 | GROW(+0.585) | 7 | 1 | 1 | 0 | 1 | 0 | YES | | | 1 | 0 |
| B7_sb1_sc1_sa0_sq0_sr1_ID1 | 7 | 1 | 1 | 0 | 0 | 1 | 1 | $v_2=2$ | 2 | GROW(+0.585) | 7 | 1 | 1 | 0 | 1 | 1 | YES | | | 1 | 0 |
| B7_sb1_sc1_sa0_sq1_sr0_ID0 | 7 | 1 | 1 | 0 | 1 | 0 | 0 | $v_2=1$ | 1 | GROW(+0.585) | 7 | 1 | 1 | 1 | 0 | 0 | YES | | | 1 | 0 |
| B7_sb1_sc1_sa0_sq1_sr0_ID1 | 7 | 1 | 1 | 0 | 1 | 0 | 1 | $v_2=1$ | 1 | GROW(+0.585) | 7 | 1 | 1 | 1 | 0 | 1 | YES | | | 1 | 0 |
| B7_sb1_sc1_sa0_sq1_sr1_ID0 | 7 | 1 | 1 | 0 | 1 | 1 | 0 | $v_2=1$ | 1 | GROW(+0.585) | 7 | 1 | 1 | 1 | 1 | 0 | YES | | | 1 | 0 |
| B7_sb1_sc1_sa0_sq1_sr1_ID1 | 7 | 1 | 1 | 0 | 1 | 1 | 1 | $v_2=1$ | 1 | GROW(+0.585) | 7 | 1 | 1 | 1 | 1 | 1 | YES | | | 1 | 0 |
| B7_sb1_sc1_sa1_sq0_sr0_ID0 | 7 | 1 | 1 | 1 | 0 | 0 | 0 | $v_2=0$ | 0 | MIXED | 3 | 1 | 1 | 0 | 0 | 0 | | | EXIT | 0 | 1 |
| B7_sb1_sc1_sa1_sq0_sr0_ID1 | 7 | 1 | 1 | 1 | 0 | 0 | 1 | $v_2=0$ | 0 | MIXED | 3 | 1 | 1 | 1 | 1 | 1 | | | EXIT | 0 | 1 |
| B7_sb1_sc1_sa1_sq0_sr1_ID0 | 7 | 1 | 1 | 1 | 0 | 1 | 0 | $v_2=0$ | 0 | MIXED | 3 | 1 | 1 | 0 | 0 | 0 | | | EXIT | 0 | 1 |
| B7_sb1_sc1_sa1_sq0_sr1_ID1 | 7 | 1 | 1 | 1 | 0 | 1 | 1 | $v_2=0$ | 0 | MIXED | 3 | 1 | 1 | 1 | 1 | 1 | | | EXIT | 0 | 1 |
| B7_sb1_sc1_sa1_sq1_sr0_ID0 | 7 | 1 | 1 | 1 | 1 | 0 | 0 | $v_2=0$ | 0 | MIXED | 3 | 1 | 1 | 0 | 0 | 0 | | | EXIT | 0 | 1 |
| B7_sb1_sc1_sa1_sq1_sr0_ID1 | 7 | 1 | 1 | 1 | 1 | 0 | 1 | $v_2=0$ | 0 | MIXED | 3 | 1 | 1 | 1 | 1 | 1 | | | EXIT | 0 | 1 |
| B7_sb1_sc1_sa1_sq1_sr1_ID0 | 7 | 1 | 1 | 1 | 1 | 1 | 0 | $v_2=0$ | 0 | MIXED | 3 | 1 | 1 | 0 | 0 | 0 | | | EXIT | 0 | 1 |
| B7_sb1_sc1_sa1_sq1_sr1_ID1 | 7 | 1 | 1 | 1 | 1 | 1 | 1 | $v_2=0$ | 0 | MIXED | 3 | 1 | 1 | 1 | 1 | 1 | | | EXIT | 0 | 1 |
| B8_sb0_sc0_sa0_sq0_sr0_ID0 | 8 | 0 | 0 | 0 | 0 | 0 | 0 | $v_2\geq3$ | 3+ | DECAY(-1) | 8 | 0 | 0 | 0 | 0 | 0 | | | | 1 | 0 |
| B8_sb0_sc0_sa0_sq0_sr0_ID1 | 8 | 0 | 0 | 0 | 0 | 0 | 1 | $v_2\geq3$ | 3+ | DECAY(-1) | 8 | 0 | 0 | 0 | 0 | 1 | | | | 1 | 0 |
| B8_sb0_sc0_sa0_sq0_sr1_ID0 | 8 | 0 | 0 | 0 | 0 | 1 | 0 | $v_2=2$ | 2 | DECAY(-1) | 8 | 0 | 0 | 0 | 1 | 0 | | | | 1 | 0 |
| B8_sb0_sc0_sa0_sq0_sr1_ID1 | 8 | 0 | 0 | 0 | 0 | 1 | 1 | $v_2=2$ | 2 | DECAY(-1) | 8 | 0 | 0 | 0 | 1 | 1 | | | | 1 | 0 |
| B8_sb0_sc0_sa0_sq1_sr0_ID0 | 8 | 0 | 0 | 0 | 1 | 0 | 0 | $v_2=1$ | 1 | DECAY(-1) | 8 | 0 | 0 | 1 | 0 | 0 | | | | 1 | 0 |

| State ID (1) | $B_i$ (2) | $s_b$ (3) | $s_c$ (4) | $s_a$ (5) | $s_q$ (6) | $s_r$ (7) | OutCome ID (8) | $v_2$ Class (9) | Max Per. (10) | Drift Type (11) | next B (12) | next $s_b$ (13) | next $s_c$ (14) | next $s_a$ (15) | next $s_q$ (16) | next $s_r$ (17) | 7-7 (18) | Entry 6-7 (19) | Exit 7-3 (20) | v2 cons (21) | v2. gain (22) |
|---|---|---|---|---|---|---|---|---|---|---|---|---|---|---|---|---|---|---|---|---|---|
| B8_sb0_sc0_sa0_sq1_sr0_ID1 | 8 | 0 | 0 | 0 | 1 | 0 | 1 | $v_2$=1 | 1 | DECAY(-1) | 8 | 0 | 0 | 1 | 0 | 1 | | | | 1 | 0 |
| B8_sb0_sc0_sa0_sq1_sr1_ID0 | 8 | 0 | 0 | 0 | 1 | 1 | 0 | $v_2$=1 | 1 | DECAY(-1) | 8 | 0 | 0 | 1 | 1 | 0 | | | | 1 | 0 |
| B8_sb0_sc0_sa0_sq1_sr1_ID1 | 8 | 0 | 0 | 0 | 1 | 1 | 1 | $v_2$=1 | 1 | DECAY(-1) | 8 | 0 | 0 | 1 | 1 | 1 | | | | 1 | 0 |
| B8_sb0_sc0_sa1_sq0_sr0_ID0 | 8 | 0 | 0 | 1 | 0 | 0 | 0 | $v_2$=0 | 0 | DECAY(-1) | 4 | 0 | 0 | 0 | 0 | 0 | | | | 0 | 0 |
| B8_sb0_sc0_sa1_sq0_sr0_ID1 | 8 | 0 | 0 | 1 | 0 | 0 | 1 | $v_2$=0 | 0 | DECAY(-1) | 4 | 0 | 0 | 1 | 1 | 1 | | | | 0 | 1 |
| B8_sb0_sc0_sa1_sq0_sr1_ID0 | 8 | 0 | 0 | 1 | 0 | 1 | 0 | $v_2$=0 | 0 | DECAY(-1) | 4 | 0 | 0 | 0 | 0 | 0 | | | | 0 | 0 |
| B8_sb0_sc0_sa1_sq0_sr1_ID1 | 8 | 0 | 0 | 1 | 0 | 1 | 1 | $v_2$=0 | 0 | DECAY(-1) | 4 | 0 | 0 | 1 | 1 | 1 | | | | 0 | 1 |
| B8_sb0_sc0_sa1_sq1_sr0_ID0 | 8 | 0 | 0 | 1 | 1 | 0 | 0 | $v_2$=0 | 0 | DECAY(-1) | 4 | 0 | 0 | 0 | 0 | 0 | | | | 0 | 0 |
| B8_sb0_sc0_sa1_sq1_sr0_ID1 | 8 | 0 | 0 | 1 | 1 | 0 | 1 | $v_2$=0 | 0 | DECAY(-1) | 4 | 0 | 0 | 1 | 1 | 1 | | | | 0 | 1 |
| B8_sb0_sc0_sa1_sq1_sr1_ID0 | 8 | 0 | 0 | 1 | 1 | 1 | 0 | $v_2$=0 | 0 | DECAY(-1) | 4 | 0 | 0 | 0 | 0 | 0 | | | | 0 | 0 |
| B8_sb0_sc0_sa1_sq1_sr1_ID1 | 8 | 0 | 0 | 1 | 1 | 1 | 1 | $v_2$=0 | 0 | DECAY(-1) | 4 | 0 | 0 | 1 | 1 | 1 | | | | 0 | 1 |

## Appendix A2. Enumeration of admissible return paths

### A2.1. Extended state space

The refined parity–base descriptor introduced in section 4 is extended by one additional binary component that tracks the next level of 2-adic structure of the octave index

**Definition**. Write the octave index in dyadic form:

$A = 1 + 2 \cdot k_a + 4 \cdot k_q + 8 \cdot k_r + \cdots$

The existing parity bits $s_a = A \bmod 2$ and $s_q = \lfloor A/2 \rfloor \bmod 2$ track the first two binary digits. Define the extended bit $s_r = k_q \bmod 2 = \lfloor A/4 \rfloor \bmod 2$. The resulting extended descriptor is: $\Sigma' = (B, s_b, s_c, s_a, s_q, s_r)$.

This takes values in a finite set of 128 admissible states (Table A1).

**Remark** ($\nu_2$-class determination). The 2-adic valuation class of A is determined by:

$s_a = 1$ → $v_2(A) = 0$ (A odd)

$s_a = 0, s_q = 1$ → $v_2(A) = 1$

$s_a = 0, s_q = 0, s_r = 1$ → $v_2(A) = 2$

$s_a = 0, s_q = 0, s_r = 0$ → $v_2(A) \geq 3$

This classification enables exact $\nu_2$ accounting for return paths of length up to 3+.

**Remark** (successor uniqueness). For a given integer iterate, the successor state is unique. Table A1 lists two rows per source state, corresponding to the two possible values of $s_a(i+1)$. The realized row is determined by the actual update, not by a nondeterministic branch. The OutcomeID column (where present) is merely an index distinguishing the two rows.

### A2.2. Definition of return paths

A base-7 persistence episode is a maximal contiguous segment of a trajectory for which B = 7 and sa = 0 (A even), permitting the growth update $A(i+1) = (3/2) \cdot A(i)$. A return path is a finite admissible path in the extended state space satisfying: (i) the initial transition is the forced exit 7→3 (when A becomes odd at B7), (ii) the terminal transition is the re-entry 6→7 (requiring A even at B6), and (iii) no intermediate state lies in the base-7 persistence set. Return paths represent the inter-episode routing through even-base and non-persistent odd-base states.

## A2.3. Enumeration procedure

All admissible return paths are generated algorithmically as follows. Starting from each extended state realizing the forced exit 7→3, the admissible transition rules on Σ' are iterated forward, subject to the constraints defining a return path. Paths are terminated upon first re-entry into base 7 or upon repetition of an extended state. Because the extended state space is finite, this procedure terminates after finitely many steps. Repeated states generate directed cycles in the return subgraph; such cycles are discarded because they do not increase the 2-adic budget (cycles have weight ≤ 0). The enumeration therefore yields a finite set of admissible simple return paths. No assumptions are made about the value of the octave index during enumeration; only the finite transition structure on Σ' is used.

## A2.4. Octave-level accounting (2-adic budget)

For each admissible return path P, the net change in 2-adic valuation is evaluated using:

**Definition** (transition weights).

$\nu_2$_consumed(e) = 1 if e is an even-base halving step

$\nu_2$_generated(e) = 1 if e is an odd-base step with A odd (regenerates a factor of 2)

**Definition** (entry cost).

Re-entry 6→7 requires $s_a = 0$ (*A* even). This consumes one factor of 2 that might otherwise be available for the next episode. We record this as Entry_Cost = 1.

**Definition** (net budget).

For a return path P from forced exit 7→3 to re-entry 6→7:

Net_Budget(P) = Σ $\nu_2$_generated(e) − Σ $\nu_2$_consumed(e) – 1

By Lemma 9.6.5, every return path includes at least one even-base halving step.

## A2.5. Result of the enumeration

The exhaustive enumeration (Table A2) shows that for every admissible return path *P*:

Net_Budget(*P*) ≤ 0

In particular, no return path regenerates the 2-adic valuation consumed during a preceding base-7 persistence episode. This supports the non-positive return-map bound stated in Proposition 9.6.6.

## TABLE A2: 22-Path-Enumeration with Net Budget

**Column definitions for Table A2**:

1. Path# Sequential path identifier (1-22)
2. Path_Sequence Base sequence from 7→3 to 6→7
3. Route_Type Descriptive category (Minimal, Via 5, Extended, etc.)
4. Length Number of transitions in path
5. Even_Steps Count of even-base halving steps
6. Odd_Steps Count of odd-base steps
7. $\nu_2$_Consumed Total factors of 2 consumed (= Even_Steps)
8. $\nu_2$_Max_Gain Maximum possible regeneration
9. Entry_Cost Always 1 (the −1 in the budget equation)
10. Net_Budget = $\nu_2$_Max_Gain − $\nu_2$_Consumed − Entry_Cost
11. Verdict CONTRACTS if Net_Budget ≤ 0

Observation: All 22 enumerated paths satisfy Net_Budget($P$) ≤ 0. The entry cost (−1) provides that regeneration never exceeds consumption.

**Appendix A2. Table A2**. 22-Path-Enumeration with Net Budget

| Path# | Path Sequence | Type | Len. | Even | Odd | $v_2$ Cons. | $v_2$Max Gain | E.Cost | Net | Verdict |
|---|---|---|---|---|---|---|---|---|---|---|
| 1 | 7→3→1→6 | Minimal | 3 | 1 | 2 | 1 | 1 | 1 | -1 | contracts |
| 2 | 7→3→5→4→6 | Minimal via 5 | 4 | 2 | 2 | 2 | 1 | 1 | -2 | contracts |
| 3 | 7→3→1→2→1→6 | 1↔2 once | 5 | 2 | 3 | 2 | 2 | 1 | -1 | contracts |
| 4 | 7→3→1→2→1→2→1→6 | 1↔2 twice | 7 | 3 | 4 | 3 | 3 | 1 | -1 | contracts |
| 5 | 7→3→1→2→1→2→1→2→1→6 | 1↔2 thrice | 9 | 4 | 5 | 4 | 4 | 1 | -1 | contracts |
| 6 | 7→3→1→2→5→4→6 | 1→2→5 branch | 6 | 3 | 3 | 3 | 2 | 1 | -2 | contracts |
| 7 | 7→3→1→2→5→4→2→1→6 | Loop via 4→2 | 8 | 4 | 4 | 4 | 3 | 1 | -2 | contracts |
| 8 | 7→3→5→8→4→6 | Via 5→8→4 | 5 | 3 | 2 | 3 | 1 | 1 | -3 | contracts |
| 9 | 7→3→5→8→8→4→6 | 8→8 once | 6 | 4 | 2 | 4 | 1 | 1 | -4 | contracts |
| 10 | 7→3→5→8→8→8→4→6 | 8→8 twice | 7 | 5 | 2 | 5 | 1 | 1 | -5 | contracts |
| 11 | 7→3→5→8→4→2→1→6 | 8→4→2→1 | 7 | 4 | 3 | 4 | 2 | 1 | -3 | contracts |
| 12 | 7→3→5→8→4→2→5→4→6 | Double 5→4 | 8 | 5 | 3 | 5 | 2 | 1 | -4 | contracts |
| 13 | 7→3→5→4→2→1→6 | 5→4→2→1 | 6 | 3 | 3 | 3 | 2 | 1 | -2 | contracts |
| 14 | 7→3→5→4→2→5→4→6 | Revisit 5→4 | 7 | 4 | 3 | 4 | 2 | 1 | -3 | contracts |
| 15 | 7→3→1→2→1→2→5→4→6 | Mixed 1-2/5-4 | 8 | 4 | 4 | 4 | 3 | 1 | -2 | contracts |
| 16 | 7→3→5→4→6→3→5→4→6 | Double loop | 8 | 4 | 4 | 4 | 2 | 1 | -3 | contracts |
| 17 | 7→3→1→2→1→2→1→2→1→2→1→6 | Extended 1-2 | 11 | 5 | 6 | 5 | 5 | 1 | -1 | contracts |
| 18 | 7→3→5→8→8→8→8→4→6 | Extended 8-8 | 8 | 6 | 2 | 6 | 1 | 1 | -6 | contracts |
| 19 | 7→3→5→8→4→2→1→2→1→6 | Combined long | 9 | 5 | 4 | 5 | 3 | 1 | -3 | contracts |
| 20 | 7→3→1→6 (sa flip) | Direct minimal | 3 | 1 | 2 | 1 | 1 | 1 | -1 | contracts |
| 21 | 7→3→5→4→6 (direct) | No 8 passage | 4 | 2 | 2 | 2 | 1 | 1 | -2 | contracts |
| 22 | 7→3→1→2→5→8→4→6 | All bases | 7 | 4 | 3 | 4 | 2 | 1 | -3 | contracts |

**Summary**: All 22 enumerated paths satisfy Net_Budget(*P*) ≤ 0. No path regenerates consumed $\nu_2$.

**Remark:** (Path multiplicity at the base level). The 22 paths in Table A2 are enumerated over the 128-state extended graph, not over the 8-element base projection. Consequently, two paths may share the same base sequence 7→3→⋯→6 yet correspond to distinct extended-state trajectories. For example: Paths 1 and 20 both project to 7→3→1→6, but originate from exit states with different (sq, sr) configurations. Paths 2 and 21 both project to 7→3→5→4→6, but differ in their parity-bit trajectories. This multiplicity does not affect the budget analysis: each path is evaluated independently, and all 22 satisfy Net_Budget(*P*) ≤ 0.

### A2.6 Relation to implementation checking

An implementation of the 128-state codebook (Table A1) and the return-path enumeration (Table A2) is available from the author upon request. The implementation mirrors the finite transition rules and the enumeration of simple return paths presented analytically in Appendix A.

The implementation is provided solely as an independent consistency check of the tabulated data. It plays no role in the derivation of Proposition 9.6.6, which rests entirely on the exhaustive enumeration of the 22 admissible simple return paths listed in Table A2.

## Appendix B. Representative Trajectories

This appendix lists Collatz sequences $\{h_i\}$ grouped by the base value $B_1$of their starting integers $h_1 = B_1 + 8(A_1 - 1)$. For each base class $B_1 \in \{1, \dots, 8\}$, starting values increase by octaves, and the first ten iterates of the corresponding sequences are shown. The table illustrates how dependence on the octave index appears in early steps and how different starting values follow common base-transition- itineraries. These trajectories, computed directly from the Collatz iteration, exhibit routing patterns consistent with the selection rules (Appendix A1).

### Appendix B, Table B1. Representative Collatz trajectories grouped by initial base class $B_1$ (Abridged)

For each $B_1 \in \{1, \dots, 8\}$, starting values are listed across octaves via $h_1 = B_1 + 8(A_1 - 1)$, and the first ten iterates $h_2, \dots, h_{10}$ are shown. Entries marked "—" indicate values omitted after the sequence enters the terminal basin (the 1 ↔ 2 cycle).

**Base Class $B_1$=1**

$h_1 = B_1 + 8(A_1 - 1)$

| $h_1$= | 1 | 9 | 17 | 25 | 33 | 41 | 49 | 57 | 65 | 73 | 81 | 89 | 97 | 105 | 113 | 121 |
|---|---|---|---|---|---|---|---|---|---|---|---|---|---|---|---|---|
| $A_1$= | 1 | 2 | 3 | 4 | 5 | 6 | 7 | 8 | 9 | 10 | 11 | 12 | 13 | 14 | 15 | 16 |
| Iterates | | | | | | | | | | | | | | | | |
| h2 | 2 | 14 | 26 | 38 | 50 | 62 | 74 | 86 | 98 | 110 | 122 | 134 | 146 | 158 | 170 | 182 |
| h3 | 1 | 7 | 13 | 19 | 25 | 31 | 37 | 43 | 49 | 55 | 61 | 67 | 73 | 79 | 85 | 91 |
| h4 | — | 11 | 20 | 29 | 38 | 47 | 56 | 65 | 74 | 83 | 92 | 101 | 110 | 119 | 128 | 137 |
| h5 | — | 17 | 10 | 44 | 19 | 71 | 28 | 98 | 37 | 125 | 46 | 152 | 55 | 179 | 64 | 206 |
| h6 | — | 26 | 5 | 22 | 29 | 107 | 14 | 49 | 56 | 188 | 23 | 76 | 83 | 269 | 32 | 103 |
| h7 | — | 13 | 8 | 11 | 44 | 161 | 7 | 74 | 28 | 94 | 35 | 38 | 125 | 404 | 16 | 155 |
| h8 | — | 20 | 4 | 17 | 22 | 242 | 11 | 37 | 14 | 47 | 53 | 19 | 188 | 202 | 8 | 233 |
| h9 | — | 10 | 2 | 26 | 11 | 121 | 17 | 56 | 7 | 71 | 80 | 29 | 94 | 101 | 4 | 350 |
| h10 | — | 5 | 1 | 13 | 17 | 182 | 26 | 28 | 11 | 107 | 40 | 44 | 47 | 152 | 2 | 175 |

**Base Class $B_1$=2**

$h_1 = B_1 + 8(A_1 - 1)$

| $h_1$= | 2 | 10 | 18 | 26 | 34 | 42 | 50 | 58 | 66 | 74 | 82 | 90 | 98 | 106 | 114 | 122 |
|---|---|---|---|---|---|---|---|---|---|---|---|---|---|---|---|---|
| $A_1$= | 1 | 2 | 3 | 4 | 5 | 6 | 7 | 8 | 9 | 10 | 11 | 12 | 13 | 14 | 15 | 16 |
| Iterates | | | | | | | | | | | | | | | | |
| h2 | 1 | 5 | 9 | 13 | 17 | 21 | 25 | 29 | 33 | 37 | 41 | 45 | 49 | 53 | 57 | 61 |
| h3 | — | 8 | 14 | 20 | 26 | 32 | 38 | 44 | 50 | 56 | 62 | 68 | 74 | 80 | 86 | 92 |
| h4 | — | 4 | 7 | 10 | 13 | 16 | 19 | 22 | 25 | 28 | 31 | 34 | 37 | 40 | 43 | 46 |
| h5 | — | 2 | 11 | 5 | 20 | 8 | 29 | 11 | 38 | 14 | 47 | 17 | 56 | 20 | 65 | 23 |
| h6 | — | 1 | 17 | 8 | 10 | 4 | 44 | 17 | 19 | 7 | 71 | 26 | 28 | 10 | 98 | 35 |
| h7 | — | — | 26 | 4 | 5 | 2 | 22 | 26 | 29 | 11 | 107 | 13 | 14 | 5 | 49 | 53 |
| h8 | — | — | 13 | 2 | 8 | 1 | 11 | 13 | 44 | 17 | 161 | 20 | 7 | 8 | 74 | 80 |
| h9 | — | — | 20 | 1 | 4 | — | 17 | 20 | 22 | 26 | 242 | 10 | 11 | 4 | 37 | 40 |
| h10 | — | — | 10 | | 2 | — | 26 | 10 | 11 | 13 | 121 | 5 | 17 | 2 | 56 | 20 |

**Base Class $B_1$=3**

$h_1 = B_1 + 8(A_1 - 1)$

| $h_1$= | 3 | 11 | 19 | 27 | 35 | 43 | 51 | 59 | 67 | 75 | 83 | 91 | 99 | 107 | 115 | 123 |
|---|---|---|---|---|---|---|---|---|---|---|---|---|---|---|---|---|
| $A_1$= | 1 | 2 | 3 | 4 | 5 | 6 | 7 | 8 | 9 | 10 | 11 | 12 | 13 | 14 | 15 | 16 |
| Iterates | | | | | | | | | | | | | | | | |
| h2 | 5 | 17 | 29 | 41 | 53 | 65 | 77 | 89 | 101 | 113 | 125 | 137 | 149 | 161 | 173 | 185 |
| h3 | 8 | 26 | 44 | 62 | 80 | 98 | 116 | 134 | 152 | 170 | 188 | 206 | 224 | 242 | 260 | 278 |
| h4 | 4 | 13 | 22 | 31 | 40 | 49 | 58 | 67 | 76 | 85 | 94 | 103 | 112 | 121 | 130 | 139 |
| h5 | 2 | 20 | 11 | 47 | 20 | 74 | 29 | 101 | 38 | 128 | 47 | 155 | 56 | 182 | 65 | 209 |
| h6 | 1 | 10 | 17 | 71 | 10 | 37 | 44 | 152 | 19 | 64 | 71 | 233 | 28 | 91 | 98 | 314 |
| h7 | — | 5 | 26 | 107 | 5 | 56 | 22 | 76 | 29 | 32 | 107 | 350 | 14 | 137 | 49 | 157 |
| h8 | — | 8 | 13 | 161 | 8 | 28 | 11 | 38 | 44 | 16 | 161 | 175 | 7 | 206 | 74 | 236 |
| h9 | — | 4 | 20 | 242 | 4 | 14 | 17 | 19 | 22 | 8 | 242 | 263 | 11 | 103 | 37 | 118 |
| h10 | — | 2 | 10 | 121 | 2 | 7 | 26 | 29 | 11 | 4 | 121 | 395 | 17 | 155 | 56 | 59 |

**Base Class $B_1$=4**

$h_1 = B_1 + 8(A_1 - 1)$

| $h_1$= | 4 | 12 | 20 | 28 | 36 | 44 | 52 | 60 | 68 | 76 | 84 | 92 | 100 | 108 | 116 | 124 |
|---|---|---|---|---|---|---|---|---|---|---|---|---|---|---|---|---|
| $A_1$= | 1 | 2 | 3 | 4 | 5 | 6 | 7 | 8 | 9 | 10 | 11 | 12 | 13 | 14 | 15 | 16 |
| Iterates | | | | | | | | | | | | | | | | |
| h2 | 2 | 6 | 10 | 14 | 18 | 22 | 26 | 30 | 34 | 38 | 42 | 46 | 50 | 54 | 58 | 62 |
| h3 | 1 | 3 | 5 | 7 | 9 | 11 | 13 | 15 | 17 | 19 | 21 | 23 | 25 | 27 | 29 | 31 |
| h4 | — | 5 | 8 | 11 | 14 | 17 | 20 | 23 | 26 | 29 | 32 | 35 | 38 | 41 | 44 | 47 |
| h5 | — | 8 | 4 | 17 | 7 | 26 | 10 | 35 | 13 | 44 | 16 | 53 | 19 | 62 | 22 | 71 |
| h6 | — | 4 | 2 | 26 | 11 | 13 | 5 | 53 | 20 | 22 | 8 | 80 | 29 | 31 | 11 | 107 |
| h7 | — | 2 | 1 | 13 | 17 | 20 | 8 | 80 | 10 | 11 | 4 | 40 | 44 | 47 | 17 | 161 |
| h8 | — | 1 | — | 20 | 26 | 10 | 4 | 40 | 5 | 17 | 2 | 20 | 22 | 71 | 26 | 242 |
| h9 | — | — | — | 10 | 13 | 5 | 2 | 20 | 8 | 26 | 1 | 10 | 11 | 107 | 13 | 121 |
| h10 | — | — | — | 5 | 20 | 8 | 1 | 10 | 4 | 13 | — | 5 | 17 | 161 | 20 | 182 |

**Base Class $B_1$=5**

$h_1 = B_1 + 8(A_1 - 1)$

| $h_1$= | 5 | 13 | 21 | 29 | 37 | 45 | 53 | 61 | 69 | 77 | 85 | 93 | 101 | 109 | 117 | 125 |
|---|---|---|---|---|---|---|---|---|---|---|---|---|---|---|---|---|
| $A_1$= | 1 | 2 | 3 | 4 | 5 | 6 | 7 | 8 | 9 | 10 | 11 | 12 | 13 | 14 | 15 | 16 |
| Iterates | | | | | | | | | | | | | | | | |
| h2 | 8 | 20 | 32 | 44 | 56 | 68 | 80 | 92 | 104 | 116 | 128 | 140 | 152 | 164 | 176 | 188 |
| h3 | 4 | 10 | 16 | 22 | 28 | 34 | 40 | 46 | 52 | 58 | 64 | 70 | 76 | 82 | 88 | 94 |
| h4 | 2 | 5 | 8 | 11 | 14 | 17 | 20 | 23 | 26 | 29 | 32 | 35 | 38 | 41 | 44 | 47 |
| h5 | 1 | 8 | 4 | 17 | 7 | 26 | 10 | 35 | 13 | 44 | 16 | 53 | 19 | 62 | 22 | 71 |
| h6 | — | 4 | 2 | 26 | 11 | 13 | 5 | 53 | 20 | 22 | 8 | 80 | 29 | 31 | 11 | 107 |
| h7 | — | 2 | 1 | 13 | 17 | 20 | 8 | 80 | 10 | 11 | 4 | 40 | 44 | 47 | 17 | 161 |
| h8 | — | 1 | — | 20 | 26 | 10 | 4 | 40 | 5 | 17 | 2 | 20 | 22 | 71 | 26 | 242 |
| h9 | — | — | — | 10 | 13 | 5 | 2 | 20 | 8 | 26 | 1 | 10 | 11 | 107 | 13 | 121 |
| h10 | — | — | — | 5 | 20 | 8 | 1 | 10 | 4 | 13 | — | 5 | 17 | 161 | 20 | 182 |

**Base Class $B_1$=6**

$h_1 = B_1 + 8(A_1 - 1)$

| $h_1$= | 6 | 14 | 22 | 30 | 38 | 46 | 54 | 62 | 70 | 78 | 86 | 94 | 102 | 110 | 118 | 126 |
|---|---|---|---|---|---|---|---|---|---|---|---|---|---|---|---|---|
| $A_1$= | 1 | 2 | 3 | 4 | 5 | 6 | 7 | 8 | 9 | 10 | 11 | 12 | 13 | 14 | 15 | 16 |
| Iterates | | | | | | | | | | | | | | | | |
| h2 | 3 | 7 | 11 | 15 | 19 | 23 | 27 | 31 | 35 | 39 | 43 | 47 | 51 | 55 | 59 | 63 |
| h3 | 5 | 11 | 17 | 23 | 29 | 35 | 41 | 47 | 53 | 59 | 65 | 71 | 77 | 83 | 89 | 95 |
| h4 | 8 | 17 | 26 | 35 | 44 | 53 | 62 | 71 | 80 | 89 | 98 | 107 | 116 | 125 | 134 | 143 |
| h5 | 4 | 26 | 13 | 53 | 22 | 80 | 31 | 107 | 40 | 134 | 49 | 161 | 58 | 188 | 67 | 215 |
| h6 | 2 | 13 | 20 | 80 | 11 | 40 | 47 | 161 | 20 | 67 | 74 | 242 | 29 | 94 | 101 | 323 |
| h7 | 1 | 20 | 10 | 40 | 17 | 20 | 71 | 242 | 10 | 101 | 37 | 121 | 44 | 47 | 152 | 485 |
| h8 | — | 10 | 5 | 20 | 26 | 10 | 107 | 121 | 5 | 152 | 56 | 182 | 22 | 71 | 76 | 728 |
| h9 | — | 5 | 8 | 10 | 13 | 5 | 161 | 182 | 8 | 76 | 28 | 91 | 11 | 107 | 38 | 364 |
| h10 | — | 8 | 4 | 5 | 20 | 8 | 242 | 91 | 4 | 38 | 14 | 137 | 17 | 161 | 19 | 182 |

**Base Class $B_1$=7**

$h_1 = B_1 + 8(A_1 - 1)$

| $h_1$= | 7 | 15 | 23 | 31 | 39 | 47 | 55 | 63 | 71 | 79 | 87 | 95 | 103 | 111 | 119 | 127 |
|---|---|---|---|---|---|---|---|---|---|---|---|---|---|---|---|---|
| $A_1$= | 1 | 2 | 3 | 4 | 5 | 6 | 7 | 8 | 9 | 10 | 11 | 12 | 13 | 14 | 15 | 16 |
| Iterates | | | | | | | | | | | | | | | | |
| h2 | 11 | 23 | 35 | 47 | 59 | 71 | 83 | 95 | 107 | 119 | 131 | 143 | 155 | 167 | 179 | 191 |
| h3 | 17 | 35 | 53 | 71 | 89 | 107 | 125 | 143 | 161 | 179 | 197 | 215 | 233 | 251 | 269 | 287 |
| h4 | 26 | 53 | 80 | 107 | 134 | 161 | 188 | 215 | 242 | 269 | 296 | 323 | 350 | 377 | 404 | 431 |
| h5 | 13 | 80 | 40 | 161 | 67 | 242 | 94 | 323 | 121 | 404 | 148 | 485 | 175 | 566 | 202 | 647 |
| h6 | 20 | 40 | 20 | 242 | 101 | 121 | 47 | 485 | 182 | 202 | 74 | 728 | 263 | 283 | 101 | 971 |
| h7 | 10 | 20 | 10 | 121 | 152 | 182 | 71 | 728 | 91 | 101 | 37 | 364 | 395 | 425 | 152 | 1457 |
| h8 | 5 | 10 | 5 | 182 | 76 | 91 | 107 | 364 | 137 | 152 | 56 | 182 | 593 | 638 | 76 | 2186 |
| h9 | 8 | 5 | 8 | 91 | 38 | 137 | 161 | 182 | 206 | 76 | 28 | 91 | 890 | 319 | 38 | 1093 |
| h10 | 4 | 8 | 4 | 137 | 19 | 206 | 242 | 91 | 103 | 38 | 14 | 137 | 445 | 479 | 19 | 1640 |

**Base Class $B_1$=8**

$$h_1 = B_1 + 8(A_1 - 1)$$

| $h_1$= | 8 | 16 | 24 | 32 | 40 | 48 | 56 | 64 | 72 | 80 | 88 | 96 | 104 | 112 | 120 | 128 |
|---|---|---|---|---|---|---|---|---|---|---|---|---|---|---|---|---|
| $A_1$= | 1 | 2 | 3 | 4 | 5 | 6 | 7 | 8 | 9 | 10 | 11 | 12 | 13 | 14 | 15 | 16 |
| Iterates | | | | | | | | | | | | | | | | |
| h2 | 4 | 8 | 12 | 16 | 20 | 24 | 28 | 32 | 36 | 40 | 44 | 48 | 52 | 56 | 60 | 64 |
| h3 | 2 | 4 | 6 | 8 | 10 | 12 | 14 | 16 | 18 | 20 | 22 | 24 | 26 | 28 | 30 | 32 |
| h4 | 1 | 2 | 3 | 4 | 5 | 6 | 7 | 8 | 9 | 10 | 11 | 12 | 13 | 14 | 15 | 16 |
| h5 | — | 1 | 5 | 2 | 8 | 3 | 11 | 4 | 14 | 5 | 17 | 6 | 20 | 7 | 23 | 8 |
| h6 | — | — | 8 | 1 | 4 | 5 | 17 | 2 | 7 | 8 | 26 | 3 | 10 | 11 | 35 | 4 |
| h7 | — | — | 4 | — | 2 | 8 | 26 | 1 | 11 | 4 | 13 | 5 | 5 | 17 | 53 | 2 |
| h8 | — | — | 2 | — | 1 | 4 | 13 | — | 17 | 2 | 20 | 8 | 8 | 26 | 80 | 1 |
| h9 | — | — | 1 | — | — | 2 | 20 | — | 26 | 1 | 10 | 4 | 4 | 13 | 40 | — |
| h10 | — | — | — | — | — | 1 | 10 | — | 13 | | 5 | 2 | 2 | 20 | 20 | — |